\documentclass[a4paper,12pt]{article}

\usepackage{amsmath}
\usepackage{amsthm}
\usepackage{amssymb}
\usepackage{cancel}
\usepackage{color}
\usepackage{calrsfs}
\definecolor{dblue}{rgb}{0.09,0.32,0.44} 
\usepackage[pdfborder={0 0 0}, colorlinks=true, pdfborderstyle={}, linkcolor=dblue, citecolor=dblue, urlcolor=blue]{hyperref}
\usepackage{epsfig}
\usepackage{psfrag}
\usepackage{subfigure}
\usepackage{graphicx}
\usepackage[mathscr,mathcal]{eucal}
\usepackage[shortlabels]{enumitem}
\usepackage{t1enc}
\usepackage[latin2]{inputenc}

\newtheorem {theorem}{Theorem}
\newtheorem {lemma}{Lemma}

\newtheorem {proposition}{Proposition}

\newtheorem* {theorem*}{Theorem}
\newtheorem* {lemma*}{Lemma}
\newtheorem* {corollary*}{Corollary}
\newtheorem* {proposition*}{Proposition}
\newtheorem* {definition*}{Definition}
\newtheorem* {conjecture*}{Conjecture}

\newtheorem* {theoremkv*} {Theorem KV}
\newtheorem* {corollarykv*} {Corollary KV}
\newtheorem* {theoremrsc*} {Theorem RSC}

\def \C {\mathbb C}

\def \N {\mathbb N}

\def \R {\mathbb R}

\def \Z {\mathbb Z}

\def\boa{\mathbf{a}}

\def\bos{\mathbf{s}}

\def\boP{\mathbf{P}}

\def\cC{\mathcal{C}}

\def\cE{\mathcal{E}}
\def\cF{\mathcal{F}}

\def\cH{\mathcal{H}}

\def\cL{\mathcal{L}}

\def\cN{\mathcal{N}}
\def\cO{\mathcal{O}}

\def\vareps{\varepsilon}

\newcommand{\probab}[1]{\ensuremath{\mathbf{P}\left(#1\right)}}
\newcommand{\expect}[1]{\ensuremath{\mathbf{E}\left(#1\right)}}
\newcommand{\var}[1]{\ensuremath{\mathbf{Var}\left(#1\right)}}

\newcommand{\condprobab}[2]{\ensuremath{\mathbf{P}\left(#1\bigm|#2\right)}}
\newcommand{\condexpect}[2]{\ensuremath{\mathbf{E}\left(#1\bigm|#2\right)}}

\newcommand{\probabom}[1]{\ensuremath{\mathbf{P}_{\omega}\left(#1\right)}}
\newcommand{\expectom}[1]{\ensuremath{\mathbf{E}_{\omega}\left(#1\right)}}
\newcommand{\varom}[1]{\ensuremath{\mathbf{Var}_{\omega}\left(#1\right)}}

\newcommand{\condprobabom}[2]{\ensuremath{\mathbf{P}_{\omega}\left(#1\bigm|#2\right)}}

\newcommand{\ind}[1]{\ensuremath{{1\!\!1}_{\{#1\}}}}

\def \toprob {\,\,\buildrel\boP\over\longrightarrow\,\,}
\def \toas {\,\,\buildrel\boa\bos\over\longrightarrow\,\,}

\DeclareMathOperator*{\infsuplim}{\overline{\underline{\lim}}}

\def \Ordo {\cO}

\def\grad {\mathrm{grad}\,}
\def\div {\mathrm{div}\,}

\def\lap{ {\mathrm{lap}}\,}

\def\Ker{\mathrm{Ker}}
\def\Ran{\mathrm{Ran}}
\def\Dom{\mathrm{Dom}}

\def\sqd{\abs{\Delta}^{1/2}}
\def\nsqd{\abs{\Delta}^{-1/2}}

\renewcommand{\d}{\mathrm d}

\newcommand{\abs}[1]{\ensuremath\left|{#1}\right|}
\newcommand{\norm}[1]{\ensuremath\left\|{#1}\right\|}
\newcommand{\sprod}[2]{\ensuremath\left\langle{#1,#2}\right\rangle}

\def \wt {\widetilde}

\def\wh{\widehat}

\title{Central Limit Theorem for Random Walks in Divergence-Free Random Drift Field: ${\cH_{-1}}$ Suffices}

\author{
{\sc Gady Kozma$^{1}$}
\qquad\qquad
{\sc B\'alint T\'oth$^{2,3}$}
\\[8pt]
{$^1$ Weizmann Institute, Rehovot, IL}
\\
{$^2$ School of Mathematics, University of Bristol, UK}
\\
{$^3$ MTA-BME Stochastics Research Group,
and}
\\
{R\'enyi Institute, Budapest, HU}
}

\begin{document}

\maketitle

\begin{abstract}
We prove central limit theorem under diffusive scaling for the displacement of a random walk on $\Z^d$ in stationary divergence-free random drift field, under the ${\cH_{-1}}$-condition imposed on the drift field. The condition is equivalent to  assuming that the stream tensor be stationary and square integrable. This improves the best existing result of Komorowski et al.\ (2012)  \cite{komorowski_landim_olla_12}, where it is assumed that the stream tensor be in $\cL^{\max\{2+\delta, d\}}$, with $\delta>0$. Our proof relies on the \emph{relaxed sector condition} of Horv\'ath et al.\ (2012) \cite{horvath_toth_veto_12}, and is technically rather simpler than existing earlier proofs of similar results by Oelschl\"ager (1988)  \cite{oelschlager_88} and Komorowski et al.\ (2012) \cite{komorowski_landim_olla_12}.



\end{abstract}

\section{Introduction: setup and main result}
\label{s:Introduction: setup and main result}

\subsection{General setup}
\label{ss:General setup}

Since its appearance in  the probability and physics literature in the mid-seventies the general topics of \emph{random walks/diffusions in random environment} became the most complex and robust area of research. For a general overview of the subject and its historical development we refer the reader to the surveys Kozlov (1985)  \cite{kozlov_85}, Zeitouni (2001)  \cite{zeitouni_01}, or Kumagai (2014)  \cite{kumagai_14}, written at various stages of this rich story. The main problem considered in our paper is that of diffusive limit in the doubly stochastic case. More precisely: we will formulate our results in a special case, that of \emph{divergence-free drift field} but emphasize that all our results and proofs go through with no extra effort, but on the expense of more notation to the general doubly stochastic case, with ellipticity condition.

Throughout this paper we denote $\cE_{d}:= \{k\in\Z^d: |k|=1\}$ the set of possible steps of a nearest-neighbour walk on $\Z^d$. A doubly stochastic, nearest-neighbour random walk in random environment is the following. Let $x\mapsto P(x) \in [0,1]^{\cE_d}$ be a stationary and ergodic (for $x\in\Z^d$) sequence of random variables, such that
\begin{align}
\label{bisto}
\sum_{k\in\cE_d} P_k(x) = 1 = \sum_{k\in\cE_d} P_{-k} (x+k),
\qquad
\mathrm{a.s.}
\end{align}
Given these,  define the random walk $X_n$ as a Markov chain on $\Z^d$, with  $X_0=0$ and conditional probabilities
\begin{align}
\label{the walk}
\condprobabom{X_{n+1}= x+k}{X_n=x} = P_k(x).
\end{align}
Denote by
\begin{align}
\label{symmetric part}
S_k(x)
&
:= \frac{P_k(x) + P_{-k}(x+k)}{2}=S_{-k}(x+k),
\\
\label{antisymmetric part}
V_k(x)
&
:= \frac{P_k(x) - P_{-k}(x+k)}{2}=V_{-k}(x+k),
\end{align}
the symmetric, respectively, antisymmetric part of the jump probabilities. Note that \eqref{bisto} is equivalent to
\begin{align}
\label{divfree 1}
\sum_{k\in\cE_d} S_k(x) \equiv 1,
\qquad
\sum_{k\in\cE_d} V_k(x) \equiv 0,
\ \ \ \
\mathrm{a.s.}
\end{align}
The second of these relations means that the lattice vector field $V$ is divergence-free, or in other words $V$ is a source- and sink-free flow on the $\Z^d$.

Assume the \emph{ellipticity} of the symmetric part: there exists a constant $\gamma>0$ so that for all $x\in \Z^d, k\in\cE_d$
\begin{align}
\label{ellipt}
S_k(x)\ge \gamma,
\qquad
\mathrm{a.s.}
\end{align}
We also assume:
\begin{align}
\label{zeromean}
\expect{V_k} =0,
\end{align}
which will ensure zero annealed mean of the walk. It is well known (see e.g. Kozlov (1985)  \cite{kozlov_85}) that due to double stochasticity \eqref{bisto} the process of \emph{the environment as seen from the position of the random walker} is stationary and ergodic in time. From \eqref{zeromean} it follows that the annealed expectation of the walk vanishes, $\expect{X_n}=\expect{\expectom{X_n}}=0$.  Under these conditions the law of large numbers
\begin{align}
\notag
\lim_{n\to\infty} n^{-1} X_n =0,
\qquad
\mathrm{a.s.}
\end{align}
follows directly from the ergodic theorem.

Regarding fluctuations around the LLN, it is not difficult to prove with bare hands that
\begin{align}
\label{diffusive lower bound}
\varliminf_{n\to \infty} n^{-1} \expect{\abs{X_n}^2} >0.
\end{align}
An upper bound on the diffusivity also follows under a subtle condition on the covariances of the drift field $\Phi: \Z^d \to \R^d(x)$,
\begin{align}
\label{Phi_def}
\Phi(x) := \sum_{k\in\cE_d} k V_k(x) = \sum_{i=1}^d e_i(V_{e_i}(x)-V_{-e_i}(x)).
\end{align}
Denoting
\begin{align}
\notag
C_{i,j}(x)
&
:=
\expect{\Phi_i(0)\Phi_j(x)},
&&
x\in\Z^d,
\\
\notag
\wh C_{i,j}(p)
&
:=
\sum_{x\in\Z^d} e^{\sqrt{-1}x\cdot p} C_{i,j}(x),
&&
p\in [-\pi,\pi)^d,
\end{align}
where the Fourier transform is meant as a distribution on $[-\pi,\pi)^d$ -- more precisely, $\wh C$ is a positive definite $d\times d$ matrix-valued measure on $[-\pi,\pi)^d$. The key condition for diffusive upper bound on the random walk $X_n$ is
\begin{align}
\label{H-1cond_0}
\int_{[-\pi,\pi)^d}
\left(\sum_{j=1}^d(1-\cos p_j) \right)^{-1} \sum_{i=1}^d \wh C_{i,i}(p) \, {\d}p <\infty.
\end{align}
Condition \eqref{H-1cond_0} is the famous $\cH_{-1}$-condition from the title of this note. Many more equivalent formulations of this condition will appear later in the paper. It is well known, see e.g. Olla (2001)  \cite{olla_01}, Komorowski et al.\ (2012)  \cite{komorowski_landim_olla_12},   that the $\cH_{-1}$-condition \eqref{H-1cond_0} implies
\begin{align}
\label{diffusive upper bound}
\varlimsup_{n\to\infty} n^{-1} \expect{\abs{X_n}^2} <\infty.
\end{align}
Now, \eqref{diffusive lower bound} and \eqref{diffusive upper bound} jointly suggest that the central limit theorem
\begin{align}
\label{clt}
n^{-1/2} X_n \Rightarrow \cN(0, \sigma^2)
\end{align}
should hold with some nondegenerate $d\times d$ covariance matrix $\sigma^2$. Attempts to prove the CLT \eqref{clt} under the (in some sense minimal) conditions \eqref{bisto}, \eqref{ellipt}, \eqref{zeromean} and \eqref{H-1cond_0} have a notorious history. In Kozlov (1985)  \cite{kozlov_85} a similar CLT is announced under the somewhat restrictive condition that the random field $x\mapsto P(x)$ be finitely dependent. As pointed out in Komorowski, Olla (2003) \cite{komorowski_olla_03a} the proof in Kozlov (1985)  \cite{kozlov_85} is incomplete. In the same paper Komorowski, Olla (2003)  \cite{komorowski_olla_03a} the CLT \eqref{clt} is stated, but as pointed out in Komorowski et al.\ (2012)  \cite{komorowski_landim_olla_12} this proof is yet again defective. Finally, in Komorowski et al.\ (2012)  \cite{komorowski_landim_olla_12} a complete proof is given, however, with more restrictive conditions: instead of the $\cH_{-1}$-condition \eqref{H-1cond_0} a rather stronger integrability condition on the field $x\mapsto P(x)$  is assumed. See the comments in section \ref{s:Historical remarks and examples}, later in this paper. More detailed historical comments on this story can be found in the notes after chapter 3 of Komorowski et al.\ (2012)  \cite{komorowski_landim_olla_12} and later in this paper. Our main result in the present paper is a complete proof of the CLT \eqref{clt}, under the conditions listed above.

The analogous diffusion problem is as follows. Let $t\mapsto X(t)\in \R^d$ be the strong solution of the SDE
\begin{align}
\label{sde}
{\d}X(t)
=
{\d}B(t) + \Phi(X(t)){\d}t,
\end{align}
where $B(t)$ is standard $d$-dimensional Brownian motion and $\Phi:\R^d\to\R^d$ is a stationary and ergodic (under space-shifts) vector field on $\R^d$ which has zero mean
\begin{align}
\notag
\expect{\Phi(x)}=0,
\end{align}
and is almost surely \emph{divergence-free}:
\begin{align}
\label{divfree in continuous space}
\div \Phi \equiv 0,
\ \ \
\mathrm{a.s.}
\end{align}
In this case the $\cH_{-1}$-condition is
\begin{align}
\label{H-1cond in continuous space}
\sum_{i=1}^d\int_{\R^d} \abs{p}^{-2}  \wh C_{i,i}(p) {\d}p <\infty,
\end{align}
where
\begin{align}
\notag
\wh C_{i,j}(p)
:=
\int_{\R^d} \expect{\Phi_i(0)\Phi_j(x)}  e^{\sqrt{-1}p\cdot x} {\d}x,
\qquad
p\in\R^d.
\end{align}
In Oelschl\"ager (1988)  \cite{oelschlager_88} it is proved that under these conditions the CLT
\begin{align}
\notag
t^{-1/2} X(t) \Rightarrow \cN(0, \sigma^2)
\end{align}
holds, with non-degenerate and finite covariance matrix $\sigma^2$. The proof is based on cutoffs and very careful technical control of approximations. It turns out that the technical details of Oelschl\"ager (1988)  \cite{oelschlager_88} cannot be transposed to the lattice walk case. For details see chapter 3 of Komorowski et al.\ (2012)  \cite{komorowski_landim_olla_12}.

Our proof, beside being complete under minimal conditions,  is also conceptually different from that of Oelschl\"ager (1988)  \cite{oelschlager_88} and
Komorowski et al.\ (2012)  \cite{komorowski_landim_olla_12} --- in our opinion more natural and technically less painful. It relies on a recent development in Kipnis-Varadhan theory of martingale approximation of additive functionals of Markov processes. We use the \emph{relaxed sector condition} of Horv\'ath et al.\ (2012) \cite{horvath_toth_veto_12} which is a natural extension of earlier sector conditions from Varadhan (1996)  \cite{varadhan_96} and
Sethuraman et al (2000)  \cite{sethuraman_varadhan_yau_00} and seems to fit well to this
problem, where we do not have a natural grading of the Hilbert space in question.

Actually, in order to keep notation simple as much as possible, we formulate our problem and result with
\begin{align}
\label{s constant}
S_k(x) \equiv \frac1{2d}.
\end{align}
It is easy to check --- and we leave this for the reader --- that (on the expense of more convolved notation) all arguments can be pushed through in the general case, assuming the ellipticity condition \eqref{ellipt}. We also note that we formulate the results in terms of continuous rather than discrete time random walk, with constant jump rate. This makes absolutely no difference in the content but lets us quote results from Kipnis-Varadhan theory in their original form.

Let us spend a few words about the physical motivation and phenomenology of the problem considered. The diffusion in divergence-free drift field, cf.  \eqref{sde}-\eqref{divfree in continuous space} may model drifting of suspended particle in stationary turbulent incompressible flow. Very similarly, the lattice counterpart \eqref{the walk} with jump rates satisfying \eqref{symmetric part}, \eqref{antisymmetric part}, \eqref{divfree 1}, \eqref{s constant} describe a random walk whose local drift is driven by a stationary source- and sink-free flow. The interest in the asymptotic description of these kind of displacement dates back to the discovery of turbulence. However, divergence-free environments appear in many other natural contexts, too. See e.g.  Komorowski et al (2012) \cite[chapter 11]{komorowski_landim_olla_12} or a surprising recent application to group theory in Bartholdi, Erschler (2011)  \cite{bartholdi_erschler_11}.

A phenomenological picture of these walks can be formulated in terms of randomly oriented cycles. Imagine that a translation invariant random ``soup of cycles'' --- that is, a Poisson point process of oriented cycles --- is placed on the lattice, and the walker is drifted along by these whirls. Now, local small cycles contribute to the diffusive behaviour. But occasionally very large cycles may cause on the long time scale faster-than-diffusive transport. Actually, this happens: in Komorowski, Olla (2002)  \cite{komorowski_olla_02} and T\'oth, Valk\'o (2012)  \cite{toth_valko_12} anomalous \emph{superdiffusive} behaviour is proved in particular cases when the $\cH_{-1}$-bound \eqref{H-1cond in continuous space} doesn't hold. Our result establishes that on the other hand, the $\cH_{-1}$-bound ensures not only boundedness of the diffusivity but also normal behaviour under diffusive scaling.

The paper is organized as follows: In the remaining subsections of the Introduction we fix notation, formulate precisely the problem, state our main result and give some comments. In section \ref{s:In the Hilbert space} we present the Hilbert space formulation, define and give some of the most important properties of the operators involved and present a form of Helmholtz's theorem which essentially says that the divergence-free drift field is the curl of the stream-tensor field. We also give an alternative formulation of the $\cH_{-1}$-condition in terms of the stream-tensor  which makes it possible to compare the content of our result with earlier ones. Section \ref{s:Historical remarks and examples} contains beside comments concrete examples. In section \ref{s:Relaxed sector condition} we recall Kipnis-Varadhan theory and the recent \emph{relaxed sector condition} which has a central role in the proof. Finally, section \ref{s:The operator B and proof of Theorem 1} contains the proof. Checking the functional-analytic condition of the abstract theorem turns into a geometric/PDE problem of having no non-trivial solution of a particular harmonic problem.

\subsection{The random walk and the \texorpdfstring{${\cH_{-1}}$}{H -1}-condition}
\label{ss:The random walk}

Let $(\Omega, \pi, \tau_z:z\in\Z^d)$ be a probability space with an ergodic $\Z^d$-action. Denote  $\cE_{d,+}:=\{e_1,\dots,e_d: e_i\in\Z^{d}, \ \ e_i\cdot e_j=\delta_{i,j}\}$ and $\cE_{d}:=\{\pm e_j: e_j\in \cE_{d,+}\} = \{k\in\Z^d: |k|=1\}$ the set of possible steps of a nearest-neighbour walk on $\Z^d$. Assume that a measurable function $v:\Omega\to\R^{\cE_{d}}$ is given with the following properties
\begin{align}
\label{vector}
&
v_k(\omega)+v_{-k}(\tau_k\omega)
\equiv0,
\\[8pt]
\label{divfree}
&
\sum_{k\in\cE_{d}}v_k(\omega)
\equiv0
\\
\label{nodrift}
&
\int_\Omega v_k(\omega) \,{\d}\pi(\omega)
=0,
\\[8pt]
\label{not ellipticity}
&
\max_{k\in\cE_d}\norm{v_k}_\infty
\le
1.
\end{align}
Define the stationary (with respect to space-shifts) lattice vector field $V:\Omega\times\Z^d \to \R^{\cE_{d}}$ by
\begin{align}
\label{vlift}
V_k(\omega, x):=v_k(\tau_x\omega).
\end{align}
Conditions \eqref{vector}, \eqref{divfree} and \eqref{not ellipticity} ensure that $V$ is almost surely a divergence-free lattice vector field on $\Z^d$ bounded by 1.  I.e.\ for all $x\in\Z^d$ and $k\in\cE_{d}$
\begin{align}
\label{lifted conditions}
V_k(x)+V_{-k}(x+k)\equiv0,
\qquad
\sum_{k\in\cE_{d}}V_k(x)\equiv0,
\qquad
\sup_{x\in\Z^d} \max_{k\in\cE_d} \abs{V_k(x)}\le1,
\end{align}
while \eqref{nodrift} implies that
\begin{align}
\notag
\expect{V_k(x)}=0,
\end{align}
which will ensure that the random walk in question has zero annealed mean.

Given $\omega\in\Omega$, define the continuous-time, nearest-neighbour Markovian random walk $t\mapsto X(t)\in \Z^d$, with $X_0=0$ and jump rates
\begin{align}
\label{jumpprob}
\condprobabom{X(t+{\d}t) = x+k }{X(t) = x}
=
\left(
1 +V_k(\omega,x)
\right) {\d}t + \Ordo(({\d}t)^2).
\end{align}
This will be a (continuous time) random walk in random environment, with divergence-free random drift field. We are primarily interested in the diffusive scaling limit of the displacement: $X(T)/\sqrt{T}$ as $T\to\infty$.

We will use the notation $\probabom{\cdot}$, $\expectom{\cdot}$  and $\varom{\cdot}$ for \emph{quenched} probability, expectation and variance. That is: probability, expectation, respectively, variance with respect to the distribution of the random walk $X(t)$, \emph{conditionally, with given fixed environment $\omega$}. The notation $\probab{\cdot}:=\int_\Omega\probabom{\cdot} {\d}\pi(\omega)$, $\expect{\cdot}:=\int_\Omega\expectom{\cdot} {\d}\pi(\omega)$ and $\var{\cdot}:=\int_\Omega\varom{\cdot} {\d}\pi(\omega) + \int_\Omega\expectom{\cdot}^2 {\d}\pi(\omega) - \expect{\cdot}^2$ will be reserved for \emph{annealed probability and expectation, etc.}. That is: probability and  expectation with respect to the random walk trajectory $X(t)$ \emph{and} the environment $\omega$.

The corresponding \emph{environment process} is the Markov process $t\mapsto \eta(t)\in \Omega$, defined by
\begin{align}
\notag
\eta(t):=\tau_{X(t)}\omega
\end{align}
This is a pure jump process on $\Omega$ with bounded (actually constant) total jump rates. So, its construction/definition does not pose any technical difficulty. It is well known (and easy to check) that due to the div-free condition \eqref{divfree} the probability measure $\pi$ is stationary and ergodic for the Markov process $t\mapsto\eta(t)$.

We will also use the variable $\varphi:\Omega\to\R^d$
\begin{align}
\label{condspeed}
\varphi(\omega)
:=
\sum_{k\in\cE_{d}} k v_k(\omega)
=
\sum_{i=1}^d e_i\left(v_{e_i}(\omega)-v_{-e_i}(\omega)\right),
\end{align}
and its stationary $\Z^d$-lifting $\Phi:\Omega\times\Z^d\to\R^d$
\begin{align}
\notag
\Phi(\omega, x):=\varphi(\tau_x\omega).
\end{align}
This is the same field as already defined in \eqref{Phi_def}. It is clear that $\expectom{{\d}X(t)}/{\d}t = \Phi(\omega, X(t))$ is the \emph{local drift} of the random walk, see \eqref{jumpprob}.

\paragraph{${\cH_{-1}}$-condition}
(first formulation): We say that $v$ satisfies the ${\cH_{-1}}$-condition if
\begin{align}
\label{H-1cond_1}
\lim_{T\to\infty} T^{-1}
\expect{\abs{\int_{0}^{T} \Phi(S(t)){\d}t}^2}
<\infty,
\end{align}
where $t\mapsto S(t)$ is a continuous-time simple symmetric random walk on $\Z^d$ with total jump rate $2d$, independent of the environment, and the expectation is  taken over the random walk \emph{and} the random scenery.

\bigskip
\noindent
Here, and further on in the paper we denote by $\abs{\dots}$ the Euclidean norm in $\R^d$. We reserve the notation $\norm{\dots}$ for the norm in the Hilbert space $\cL^2(\Omega,\pi)$ to be introduced later.

As discussed in the previous subsection, there is an equivalent formulation of the ${\cH_{-1}}$-condition is in terms of the covariance matrix of the random field $\Phi(\omega,x)$. Define $C:\Z^d\to\R^{d\times d}$
\begin{align}
\label{covariance matrix of Phi}
C_{i,j}(x)
:=
&
\expect{\Phi_i(x)\Phi_j(0)}
\end{align}
Note that $\Z^d\ni x\mapsto C(x)\in \R^{d \times d}$ is of positive type: for any $n\in\N$, $x_1,\dots,x_n\in\Z^d$, $y_1,\dots,y_n\in\C$ and $z_1,\dotsc,z_d\in\C$,
\begin{align}
\notag
\sum_{k,l=1}^n \sum_{i,j=1}^d \overline{y_kz_i} C_{i,j}(x_k-x_l) y_lz_j\ge0.
\end{align}
Thus, it follows from Bochner's Theorem (see Theorem IX.9 in Reed, Simon (1975)  \cite{reed_simon_vol1_vol2_75}) that the Fourier transform of $\sum \overline{z_i}C_{i,j}z_j$ is a positive measure for all $z$, and hence $\wh C$ is a ($d\times d$)-positive-definite-matrix-valued measure on $[-\pi,\pi]^d$. With some abuse of notation we will denote
\begin{align}
\label{Fourier transform of covariance matrix of Phi}
\wh C_{i,j}(p)
&:=
\sum_{x\in\Z^d} e^{\sqrt{-1}p\cdot x} C_{i,j}(x).
\end{align}
That is, we denote by $\wh C_{i,j}(p)\, \d{}p$ the Fourier transform of  the covariance matrix $C_{i,j}(x)$, although it may not be absolutely continuous with respect to the Lebesgue measure ${\d}p$ on $[-\pi,\pi]^d$.

Let $\wh D:[-\pi,\pi]^d\to \R$ be the Fourier transform of the $\Z^d$-Laplacian
\begin{align}
\notag
\wh D(p):=\sum_{i=1}^d (1-\cos(p\cdot e_i)),
\end{align}
and
\begin{align}
\label{covariance of RWRS}
\wt C_{i,j}
:=
\frac{2}{(2\pi)^{d}}
\int_{[-\pi,\pi]^d}
\wh D(p)^{-1}
\wh C_{i,j}(p) {\,\d}p,
\qquad
i,j=1,\dots, d.
\end{align}
It is easy to see that $\wt C_{i,j}$ is the asymptotic covariance matrix of $T^{-1/2}\int_0^T \Phi(S(t)) {\,\d}t \in\R^d$. Indeed,
\begin{align}
\notag
\lefteqn{
\lim_{T\to\infty} T^{-1}
\expect{ \int_{0}^{T} \Phi_i(S(t)){\,\d}t \int_{0}^{T} \Phi_j(S(t)){\,\d}t}
}\hskip3cm
&
\\
\notag
\textrm{By stationarity}\qquad&
=
\lim_{T\to\infty}
\int_{0}^{T}
\frac{T-t}{T}
\left(
\expect{
\Phi_i(0)
\Phi_j(S(t))
}
+
\expect{
\Phi_j(0)
\Phi_i(S(t))
}
\right)
{\,\d}t
\\
\notag
&
=
2\lim_{T\to\infty}
\int_{0}^{T}
\frac{T-t}{T}
\sum_{x\in\Z^d}
\probab{S_t=x} C_{i,j}(x)
{\,\d}t
\\
\notag
\textrm{By Parseval}\qquad&
=
\frac{2}{(2\pi)^d}
\lim_{T\to\infty}
\int_{0}^{T}
\frac{T-t}{T}
\int_{[-\pi,\pi]^d}
e^{-t\wh D(p)} \wh C_{i,j}(p)\, {\d}p\, {\d}t
\\
\notag
\textrm{By Fubini \& Fatou}\qquad
&=\frac{2}{(2\pi)^d}
\int_{[-\pi,\pi]^d} \wh C_{i,j}(p)
\lim_{T\to\infty}
\int_0^T\frac{T-t}{T}e^{-t\wh D(p)} \, {\d}t\, {\d}p
\\
\label{covariance_matrix of RWRS}
\textrm{By Fubini}\qquad
&
=
\frac{2}{(2\pi)^d}
\int_{[-\pi,\pi]^d}
\wh C_{i,j}(p)\wh D(p)^{-1} \, {\d}p
=
\wt C_{i,j}
\end{align}
where the applications of Fubini and Fatou's theorems are justified by positivity and monotonicity. Hence the alternative equivalent formulation of the ${\cH_{-1}}$-condition:

\paragraph{${\cH_{-1}}$-condition}
(second formulation):
\begin{align}
\label{H-1cond_2}
(2\pi)^{-d}
\int_{[-\pi,\pi]^d}
\sum_{i=1}^d
\wh D(p)^{-1}
\wh C_{i,i}(p) {\d}p <\infty.
\end{align}

\bigskip
\noindent
Note, that this is the same as \eqref{H-1cond_0}. We will see below \emph{yet another} equivalent formulation of the ${\cH_{-1}}$-condition, \eqref{H-1cond_3}, in functional analytic terms,  in the spirit of Kipnis-Varadhan theory. We will refer to either one of the equivalent formulations as the \emph{${\cH_{-1}}$-condition} for the drift field.

\subsection{Central limit theorem for the random walk}
\label{ss:Central limit theorem for the random walk}

As already mentioned, we are primarily interested in the asymptotic diffusive behaviour of the walk $t\mapsto X(t)$. By the definition of the drift $\varphi$ (recall \eqref{condspeed}),
\begin{align}
\label{forward}
\condexpect{X(t+dt)-X(t)}{\eta(t)=\omega}
=
\varphi(\omega)dt + \Ordo((dt)^2)
\end{align}
and divergence freeness gives
\begin{align}
\label{backward}
\condexpect{X(t-dt)-X(t)}{\eta(t)=\omega}
=
-
\varphi(\omega)dt + \Ordo((dt)^2),
\end{align}
with $\varphi:\Omega\to\R^d$ the drift defined in \eqref{condspeed}. So, it is most natural to decompose $X(t)$ as
\begin{align}
\notag
X(t)=Y(t)+Z(t),
\end{align}
where
\begin{align}
\label{Y and Z def}
&
Y(t)
:=
X(t)-
\int_{0}^{t} \varphi(\eta(s)){\d}s,
&&
Z(t)
:=
\int_{0}^{t} \varphi(\eta(s)){\d}s.
\end{align}
Due to \eqref{forward} $t\mapsto Y(t)$ is a martingale.

The following diffusive (lower and upper) bounds are easy consequences of this decomposition:

\begin{proposition}
\label{prob:diffusive}
Assume that the random drift $\varphi$ satisfies the ${\cH_{-1}}$-condition \eqref{H-1cond_1}/\eqref{H-1cond_2}. Then the following diffusive lower and upper bounds hold for the random walk $t\mapsto X(t)$ defined in \eqref{jumpprob}:
\begin{align}
\label{diffusive bounds}
2d
\le
\infsuplim_{T\to\infty} T^{-1} \expect{\abs{X(T)}^2}
\le
2d
+
8\sum_{i=1}^d \wt C_{i,i}
<
\infty
\end{align}
where $\wt C$ is defined by \eqref{covariance of RWRS}.
\end{proposition}

\begin{proof}
Due to \eqref{forward} and \eqref{backward} the processes $t\mapsto Y(t)$ and $t\mapsto Z(t)$ are uncorrelated under the annealed measure. Hence
\begin{align}
\notag
\expect{X_i(t)X_j(t)}
=
\expect{Y_i(t)Y_j(t)}
+
\expect{Z_i(t)Z_j(t)}.
\end{align}
Because the rates $1+V$ sum to $2d$, $\sum_{i=1}^d\expect{Y_i(t)^2}=2dt$, hence the lower bound in \eqref{diffusive bounds}.
On the other hand, it is well known, see e.g. Komorowski et al.\ (2012)  \cite{komorowski_landim_olla_12}, or  Olla (2001) \cite[proposition (2.1.11)]{olla_01}, that
\begin{align}
\label{variance bound on Z}
\expect{\abs{Z(T)}^2}
=
\expect{\abs{\int_{0}^{T} \Phi(X(t)){\d}t}^2}
\le
8
\expect{\abs{\int_{0}^{T} \Phi(S(t)){\d}t}^2}.
\end{align}
Dividing by $T$ and taking the limit $T\to\infty$ in \eqref{variance bound on Z} we obtain --- via \eqref{covariance_matrix of RWRS} --- the upper bound in \eqref{diffusive bounds}.
\end{proof}

Actually, it is not just the upper bound \eqref{variance bound on Z} that holds but the covariance matrix $\expect{Z_i(t)Z_j(t)}$ is dominated (as $d\times d$ matrix) by the matrix $8t \wt C_{i,j}$ defined in \eqref{covariance of RWRS}.

The main result of the present paper is the following:

\begin{theorem}
\label{thm:main}
If the \emph{${\cH_{-1}}$-condition} \eqref{H-1cond_2} holds then the asymptotic covariance matrix
\begin{align}
\notag
(\sigma^2)_{i,j}
:=
\lim_{T\to\infty} T^{-1}\expect{X_i(T)X_j(T)}
\end{align}
exists,
\begin{align}
\label{bounds on sigmasquared}
2d I_d
\le
\sigma^2
\le
2d I_d
+
8 \wt C,
\end{align}
and for any $m\in\N$, $t_1,\dots,t_m\in\R_+$ and any test function $F:\R^{md}\to\R$ continuous and bounded
\begin{align}
\label{eq:quannealed}
\lim_{T\to\infty}
\int_\Omega
\abs{
\expectom{F\left(\frac{X(Tt_1)}{\sqrt{T}}, \dots, \frac{X(Tt_m)}{\sqrt{T}}\right)}
-
\expect{F(W(t_1), \dots, W(t_m))}
} {\d}\pi(\omega)=0,
\end{align}
where $t\mapsto W(t)\in\R^d$ is a Brownian motion with
\begin{align}
\notag
&
\expect{W_i(t)}=0,
&&
\expect{W_i(s)W_j(t)}=\min\{s,t\}(\sigma^2)_{i,j}
\end{align}
\end{theorem}

\subsubsection*{Remarks:}

\begin{enumerate}[(1)]

\item
Since the sum of the jump rates on the right hand side of \eqref{jumpprob} is constant $2d$, there is no difference whatsoever between discrete or continuous time. We choose the continuous time formulation only for reasons of convenience: some formulas become somewhat simpler and the \emph{relaxed sector condition} on which our proof relies was in Horv\'ath et al.\  (2012) \cite{horvath_toth_veto_12} originally formulated in the (a priori more general) continuous time setting.

\item
Note that the inequality in \eqref{not ellipticity} is \emph{not strict}: this is \emph{not an ellipticity condition}. It is allowed that jumps along some directed edges be prohibited (and the jump rates along the same edge oppositely directed be maximal, 2).

\item
\label{general doubly stochastic}
Let us stress again that the most general case of nearest-neighbour continuous-time random walk in \emph{doubly stochastic random environment} is
\begin{align}
\label{jumpprob_gen}
\condprobabom{X(t+dt) = x+k }{X(t) = x}
=
p_k(\tau_x\omega) dt + \Ordo((dt)^2),
\end{align}
with
\begin{align}
\label{bistochastic}
\sum_{k\in\cE_d} p_k(\omega)
=
\sum_{k\in\cE_d} p_{-k}(\tau_k\omega).
\end{align}
In this general case it is natural to write $p_k(\omega)=s_k(\omega)+v_k(\omega)$, with
\begin{align}
\notag
&
s_k(\omega)
:=
\frac{p_k(\omega)+p_{-k}(\tau_k\omega)}{2},
&&
v_k(\omega)
:=
\frac{p_k(\omega)-p_{-k}(\tau_k\omega)}{2}.
\end{align}
In the present note we analyse the case when $s_k(\omega)\equiv1$, for all $k\in\cE_d$. However, with the expanse of more notation and length of exposition we can do exactly the same in the general case, under the \emph{ellipticity} condition
\begin{align}
\label{ellipticity}
s_k(\omega)
\ge
\epsilon
>0,
\qquad
\forall \ k\in\cE_d.
\end{align}

\item
Note that the theorem establishes a CLT for the displacement of the
random walker ``halfway'' between annealed and quenched. Annealed CLT
would have been
\begin{align}
\notag
\lim_{T\to\infty}
&
\abs{\expect{F \left( \frac{X(tT)}{\sqrt{T}} \right)}- \expect{F(W(t))}}
=
\\
\notag
&\hskip2cm
\lim_{T\to\infty}\abs{\int_\Omega \expectom{F\left( \frac{X(tT)}{\sqrt{T}} \right) }\,\d\pi(\omega)-\expect{F(W(t)) }}
=0,
\end{align}
which is weaker than what we get in \eqref{eq:quannealed}, while a quenched CLT
would be
\begin{align}
\notag
\lim_{T\to\infty}
\expectom{ F\left( \frac{X(tT)}{\sqrt{T}} \right)}
=
\expect{F(W\dotsc)},
\qquad \pi-\text{a.s.},
\end{align}
which is stronger than \eqref{eq:quannealed}.
The limit in \eqref{eq:quannealed} is the type of result which drops out from the martingale approximation a la Kipnis-Varadhan. See e.g.\ definition (2.2.a) and Theorem 1 in T\'oth (1986)  \cite{toth_86}. See also Corollary KV in section \ref{s:Relaxed sector condition} of this paper.

\item
The analogous model in continuous space-time  is the diffusion $t\mapsto X(t)\in\R^d$ in an ergodic (with respect to spacial translations), almost surely divergence free and zero mean drift field $\Phi:\R^d\to\R^d$, as specified in \eqref{sde}, \eqref{divfree in continuous space}. Analogous result to Theorem \ref{thm:main} can be proved in a very similar way -- reproducing this way a conceptually and technically simpler proof of Theorem 1 of Oelschl\"ager (1988) \cite{oelschlager_88}. However, in order to keep this paper reasonably short we will not give here details for the diffusion model \eqref{sde}, \eqref{divfree in continuous space}.

\item
More comments on the place of Theorem \ref{thm:main} and its relation to other results will be formulated later in section \ref{s:Historical remarks and examples}, after a necessary discussion of the stream tensor.

\item
Theorem \ref{thm:main} is a step towards the following conjecture,
inspired by the work of Varadhan: let $X$ be a random walk in a
stationary, ergodic, elliptic random environment. Then if $X$ is
diffusive (in the sense that its variance grows linearly with time),
then it also satisfies the central limit theorem. The $\cH_{-1}$
condition is equivalent to diffusivity in the \emph{reversible} case,
but not in general: Peled (private communication) gave an example of
a diffusive divergence-free random walk which is not in $\cH_{-1}$.

\end{enumerate}

\section{In the Hilbert space \texorpdfstring{$\cL^2(\Omega,\pi)$}{}}
\label{s:In the Hilbert space}

\subsection{Spaces and operators}
\label{ss:Spaces and operators}

It is most natural to put ourselves into the Hilbert space
\begin{align}
\notag
\cH:=\Big\{f\in\cL^2(\Omega,\pi):\int_\Omega f {\d}\pi=0\Big\}.
\end{align}
We denote by $T_x$, $x\in\Z^d$, the spatial shift operators
\begin{align}
\notag
T_x f(\omega)
:=
f(\tau_x\omega),
\end{align}
and note that they are unitary:
\begin{align}
\notag
T_x^*=T_{-x}=T_x^{-1}.
\end{align}
The $\cL^2$-gradients $\nabla_k$, $k\in\cE_{d}$, respectively, $\cL^2$-Laplacian $\Delta$, are:
\begin{align}
\notag
&
\nabla_k
:=
T_k-I,
&&
\nabla_k^*
=
\nabla_{-k},
&&
\norm{\nabla_k}=2,
\\[8pt]
\notag
&
\Delta
:=
\sum_{l\in\cE_{d}}\nabla_l
=
-
\frac{1}{2}
\sum_{l\in\cE_{d}}\nabla_l\nabla_{-l},
&&
\Delta^*
=
\Delta\le0,
&&
\norm{\Delta}=4d.
\end{align}
Due to ergodicity of the $\Z^d$-action $(\Omega, \pi, \tau_z:z\in\Z^d)$,
\begin{align}
\label{kerDelta}
\Ker(\Delta)=\{0\}.
\end{align}
Indeed, $\Delta f=0$ implies that $0=\sprod{f}{\Delta f} = \sum_{i=1}^d \sprod{\nabla_i f}{\nabla_i f}$ and since all terms are non-negative, they must all be 0 and $f$ must be invariant to translations. Ergodicity to $\Z^d$ actions means that $f$ is constant, and since our Hilbert space is that of functions averaging to zero, $f$ must be zero.

We define the bounded positive operator $\sqd$ in terms of the spectral theorem (applied to the bounded positive operator $\abs{\Delta}:=-\Delta$). Note that due to \eqref{kerDelta} $\Ran \abs{\Delta}$ is dense in $\cH$, and hence so is $\Ran\abs{\Delta}^{1/2}$ which is a superset of it. Hence it follows that $\nsqd :=\left(\sqd\right)^{-1}$ is an (unbounded) positive self-adjoint operator with $\Dom \nsqd = \Ran \sqd$ and $\Ran \nsqd = \Dom \sqd = \cH$. Note that the dense subspace $\Dom \nsqd = \Ran \sqd$ is invariant under, and the operators $\sqd$ and $\nsqd$ commute with the  translations $T_x$, $x\in\Z^d$.

We define, for all $k\in\cE$, the \emph{Riesz operators}
\begin{align}
\label{Gamma_k}
\Gamma_k:\Dom \nsqd \to \cH,
\qquad
\Gamma_k=\Delta^{-1/2}\nabla_k=\nabla_k\Delta^{-1/2},
\end{align}
and note that for any $f\in\Dom \nsqd $
\begin{align}
\notag
\norm{\Gamma_k f}^2
=
\sprod{\nsqd f}{\nabla_{-k}\nabla_k \nsqd f}
\le
\sprod{\nsqd f}{\abs {\Delta} \nsqd f}
=
\norm{f}^2.
\end{align}
Thus, the operators $\Gamma_k$, $k\in\cE$, extend as bounded operators to the whole space $\cH$. The following properties are easy to check:
\begin{align}
\label{Gammaadj}
&
\Gamma_k^*=\Gamma_{-k},
&&
\norm{\Gamma_k}=1,
&&
\frac12
\sum_{l\in\cE_{d}}\Gamma_l\Gamma_l^*=I.
\end{align}

A third equivalent formulation of the ${\cH_{-1}}$-condition is the following:

\paragraph{${\cH_{-1}}$-condition}
(third formulation):
For all $k\in\cE$,
\begin{align}
\label{H-1cond_3}
v_k\in\Dom \nsqd = \Ran \sqd.
\end{align}

\bigskip
\noindent
To see the equivalence first note that our first version of the $\cH_{-1}$ condition, \eqref{H-1cond_1} is equivalent to $\varphi_i\in\Dom \nsqd$ for $i=1,\dots,d$. Indeed, in operator theory language the distribution of $\Phi_i(S(t))$ is the lifting to $\Z^d$ of $e^{-t\Delta}\varphi_i$ and hence a calculation similar to (\ref{covariance_matrix of RWRS}) shows that
\begin{align}
\notag
\frac 1T
\expect{\abs{\int_0^T\Phi(S(t))\,{\d}t}^2}
=
\sum_{i=1}^d
\sprod{\varphi_i}{\int_0^T\frac{T-t}{T}e^{-t\Delta}\varphi_i}
\end{align}
and an application of the spectral theorem for $|\Delta|$ shows that this is bounded in $T$ if and only if all $\varphi_i$ are in $\Dom|\Delta|^{-1/2}$. To conclude from $\varphi$ to $v$ we recall that $\varphi_i = (I+T_{e_i})v_{e_i}=(2I+\nabla_{e_i})v_{e_i}$. Since $\Gamma_{e_i}=|\Delta|^{-1/2}\nabla_{e_i}$ is bounded,  we get that $\nabla_{e_i} v_{e_i} \in \Dom(|\Delta|^{-1/2})$. Rearranging gives
\begin{align}
\notag
\varphi_i-2v_i \in \Dom(|\Delta|^{-1/2})
\end{align}
which shows that $\varphi_i\in\Dom(|\Delta|^{-1/2})$ if and only if so is $v_i$.

Finally, we also define the multiplication operators $M_k$, $k\in\cE_d$,
\begin{align}
\label{vmultipl}
&
M_k f(\omega):= v_k(\omega) f(\omega),
&&
M_k^*=M_k,
&&
\norm{M_k}\le1.
\end{align}
The following commutation relation holds due to the div-free condition \eqref{divfree}, and is easily checked:
\begin{align}
\label{Mnablacommute}
\sum_{l\in\cE_{d}} M_l \nabla_l + \sum_{l\in\cE_{d}} \nabla_{-l} M_{l} =0.
\end{align}

The \emph{infinitesimal generator} of the stationary environment process $t\mapsto\eta(t)$, acting on the Hilbert space $\cL^2(\Omega,\pi)$ is:
\begin{align}
\notag
L= \Delta + \sum_{l\in\cE_{d}} M_l \nabla_l.
\end{align}
We decompose the infinitesimal generator $L$ as sum of its self-adjoint and skew self-adjoint part, $L=-S+A$. Due to the commutation relation \eqref{Mnablacommute}:
\begin{align}
\label{S and A}
&
-S
:=
(L+L^*)/2
=
\Delta,
&&
A
:=
(L-L^*)/2
=
\sum_{l\in\cE_{d}} M_l \nabla_l
=
-\sum_{l\in\cE_{d}} \nabla_{-l} M_{l}.
\end{align}

\subsection{The stream tensor field}
\label{ss:The stream tensor field}

The following proposition establishes the existence of the stream
tensor field and is essentially Helmholtz's theorem. It is the $\Z^d$
lattice counterpart of Proposition 11.1 from
Komorowski et al.\ (2012)  \cite{komorowski_landim_olla_12}. It is actually not necessary in a technical sense for the proof of
our main result, but it is an important part of the story and sheds light on the role and limitations of the ${\cH_{-1}}$-condition in this context.

\begin{proposition}
\label{prop:helmholtz}
(i)
There exists an antisymmetric tensor field
${H}:\Omega\times\Z^d\to\R^{\cE_{d}\times\cE_{d}}$ such that
for all $x\in\Z^d$ we have ${H}_{k,l}(\cdot,x)\in\cH$ and
\begin{align}
\label{Thetasymms}
{H}_{l,k}(\omega,x)
=
{H}_{-k,l}(\omega,x+k)
=
{H}_{k,-l}(\omega,x+l)
=
-{H}_{k,l}(\omega,x),
\end{align}
with stationary increments
\begin{align}
\notag
{H}(\omega,y)-{H}(\omega,x)={H}(\tau_x\omega, y-x)-{H}(\tau_x\omega, 0),
\end{align}
such that
\begin{align}
\label{V=curlTheta}
V_k(\omega,x)=\sum_{l\in\cE_{d}} {H}_{k,l}(\omega,x).
\end{align}
The realization of the tensor field ${H}$ is uniquely determined by the ``pinning down'' condition \eqref{Theta0} below.

\smallskip
\noindent
(ii)
The ${\cH_{-1}}$-condition \eqref{H-1cond_3} holds if and only if there exist ${h}_{k,l}\in\cH$, $k,l\in\cE_d$,  such that
\begin{align}
\label{thetatensor}
{h}_{l,k}
=
T_k {h}_{-k,l}
=
T_l {h}_{k,-l}
=
-
{h}_{k,l}
\end{align}
and
\begin{align}
\label{v=curltheta}
v_k(\omega)=\sum_{l\in\cE_{d}}{h}_{k,l}(\omega).
\end{align}
In this case the tensor field ${H}$ can be realized as the stationary  lifting of ${h}$:
\begin{align}
\label{lifting of theta}
{H}_{k,l}(\omega,x)={h}_{k,l}(\tau_x\omega).
\end{align}
\end{proposition}

\begin{proof}
(i)
For $k,l,m\in\cE_{d}$ define
\begin{align}
\notag
g_{m;k,l}
:=
\Gamma_m\big(\Gamma_l v_k - \Gamma_k v_l\big),
\end{align}
where $\Gamma_l$ are the Riesz operators defined in \eqref{Gamma_k}, and note that for all $k,l,m,n\in\cE_{d}$
\begin{align}
\label{gistensor}
&
g_{m;l,k}
=
T_kg_{m;-k,l}
=
T_lg_{m;k,-l}
=
-g_{m;k,l},
\\[8pt]
\label{gisgrad}
&
g_{m;l,k}+T_mg_{n;l,k}
=
g_{n;l,k}+T_ng_{m;l,k},
\\[8pt]
\label{curlg=gradv}
&
\sum_{l\in\cE_{d}}g_{m;k,l} = \nabla_m v_k.
\end{align}
\eqref{gistensor} means that that keeping the subscript $m\in\cE_{d}$ fixed, $g_{m;k,l}$ has exactly the symmetries of a $\cL^2$-tensor variable indexed by  $k,l\in\cE_{d}$. \eqref{gisgrad} means that, on the other hand, keeping $k,l\in\cE_{d}$ fixed, $g_{m;k,l}$ is a $\cL^2$-gradient in the subscript $m\in\cE_{d}$. Finally, \eqref{curlg=gradv} means that the $\cL^2$-curl of tensor $g_{m;\cdot,\cdot}$ is actually the $\cL^2$-gradient of the vector $v_\cdot$.

Let  $G_{m;k,l}:\Omega\times\Z^d\to\R$ be the lifting  $G_{m;k,l}(\omega,x):=g_{m;k,l}(\tau_x\omega)$. By \eqref{gisgrad}, for any $k,l\in\cE_{d}$ fixed $\left(G_{m;k,l}(\omega,x)\right))_{m\in\cE_{d}}$ is a lattice gradient. The increments of ${H}_{k,l}$ are defined by
\begin{align}
\label{gradTheta=G}
{H}_{k,l}(\omega,x+m)-{H}_{k,l}(\omega,x)
=
G_{m;k,l}(\omega,x),
\ \ \ m\in\cE_{d}.
\end{align}
This is consistent, due to \eqref{gisgrad}.

Next, in order to uniquely determine the tensor field $H$,  we ``pin down'' its values at $x=0$.
For $e_i,e_j\in\cE_{d,+}$ choose
\begin{align}
\label{Theta0}
\begin{aligned}
&
{H}_{e_i,e_j}(\omega,0)
=
0,
&&
{H}_{-e_i,e_j}(\omega,0)
=
-g_{-e_i;e_i,e_j}(\omega),
\\[8pt]
&
{H}_{e_i,-e_j}(\omega,0)
=
g_{-e_j;e_i,e_j}(\omega),
&&
{H}_{-e_i,-e_j}(\omega,0)
=
-g_{-e_i;e_i,e_j}(\omega)+g_{-e_j;e_i,e_j}(\tau_{-e_i}\omega).
\end{aligned}
\end{align}
The tensor field ${H}$ is fully determined by \eqref{gradTheta=G} and \eqref{Theta0}. Due to \eqref{gistensor} and \eqref{curlg=gradv}, \eqref{Thetasymms}, respectively, \eqref{V=curlTheta} will hold, indeed.

(ii)
First we prove the only if part. Assume \eqref{H-1cond_3} and let
\begin{align}
\notag
{h}_{k,l}
=
\Gamma_l \abs{\Delta}^{-1/2} v_k - \Gamma_k  \abs{\Delta}^{-1/2} v_l
=
\abs{\Delta}^{-1/2} \big( \Gamma_l  v_k - \Gamma_k  v_l\big).
\end{align}
Hence \eqref{thetatensor} and \eqref{v=curltheta} are readily obtained. Next we prove the if part. Assume that there exist $h_{k,l}\in\cH$ with the symmetries \eqref{thetatensor} and $v_k$ is realized as in \eqref{v=curltheta}. Then we have
\begin{align}
\notag
v_k
=
\sum_{l\in\cE_d} h_{k,l}
=
\frac12 \sum_{l\in\cE_d} \left(h_{k,l} + h_{k,-l}\right)
=
-\frac12 \sum_{l\in\cE_d} \nabla_l h_{k,-l}
=
-\frac12 \abs{\Delta}^{1/2} \sum_{l\in\cE_d} \Gamma_l h_{k,-l},
\end{align}
which shows indeed \eqref{H-1cond_3}.
\end{proof}

\paragraph{${\cH_{-1}}$-condition}
(fourth formulation):
The drift vector field $V$ is realized as the curl of a \emph{stationary and square integrable}, zero mean  tensor field ${H}$, as shown in \eqref{V=curlTheta}.

\bigskip
\noindent
{\bf Remark.}
If the $\cH_{-1}$-condition \eqref{H-1cond_1}/\eqref{H-1cond_2}/\eqref{H-1cond_3} does not hold it may still be possible that there exists a \emph{non-square integrable} tensor variable ${h}:\Omega\to\R^{\cE_{d}\times\cE_{d}}$ which has the symmetries \eqref{thetatensor} and with $v:\Omega\to\R^{\cE_{d}}$ realized as in \eqref{v=curltheta}. Then let ${H}:\Omega\times\Z^d\to\R^{\cE_{d}\times\cE_{d}}$ be the stationary lifting \eqref{lifting of theta} and we still get \eqref{V=curlTheta} with a \emph{stationary but not square integrable} tensor field. Note that this is not decidable in terms of the covariance matrix \eqref{covariance matrix of Phi} or its Fourier transform \eqref{Fourier transform of covariance matrix of Phi}. The question of diffusive (or super-diffusive) asymptotic behaviour of the walk $t\mapsto X(t)$ in these cases is fully open.

\medskip

In the next proposition -- which essentially follows an argument from the otherwise incomplete proof of Theorem II.3.3 in Kozlov (1985)  \cite{kozlov_85} -- we give a sufficient condition for the ${\cH_{-1}}$-condition \eqref{H-1cond_2} to hold.

\begin{proposition}
\label{prop:suff cond for H-1}
If $p\mapsto \wh C(p)$ is twice continuously differentiable function in a neighbourhood of $p=0$ then the ${\cH_{-1}}$-condition \eqref{H-1cond_2} holds.
\end{proposition}

\begin{proof}
For the duration of this proof we introduce the notation
\begin{align}
\notag
B_{k,l}(x):=\expect{V_k(0)V_l(x)},
\quad
\wh B_{k,l}(p) := \sum_{x\in \Z^d} e^{\sqrt{-1}x\cdot p} B_{k,l} (x),
\end{align}
with $k,l\in\cE_d, x\in\Z^d, p\in[-\pi,\pi]^d$. Hence (recall the definition of $C_{i,j}$,
\eqref{Phi_def}) for $i,j\in\{1,\dots,d\}$
\begin{align}
\notag
\wh C_{i,j}(p)
=
\wh B_{e_i,e_j}(p) - \wh B_{-e_i,e_j}(p) - \wh B_{e_i,-e_j}(p) + \wh B_{-e_i,-e_j}(p).
\end{align}
(The identity is meant in the sense of distributions.)

Note that  due to the first clause in \eqref{lifted conditions}
\begin{align}
\label{Chatvector}
\wh B_{k,l}(p)
=
-e^{\sqrt{-1}p\cdot k}\wh B_{-k,l}(p)
=
-e^{-\sqrt{-1}p\cdot l}\wh B_{k,-l}(p)
=
e^{\sqrt{-1}p\cdot (k-l)}\wh B_{-k,-l}(p).
\end{align}
Using \eqref{Chatvector} in the above expression of $C(p)$ in terms of $B(p)$, direct computations yield
\begin{align}
\notag
\wh C_{i,j} = \left(1+e^{-\sqrt{-1}p\cdot e_i}\right) \left(1+e^{\sqrt{-1}p\cdot e_j}\right) \wh B_{e_i,e_j}(p).
\end{align}
Thus, the regularity condition imposed on $p\mapsto C(p)$ is equivalent to assuming the same regularity about $p\mapsto \wh B(p)$.

Next, due to the second clause of \eqref{lifted conditions}
\begin{align}
\label{Chatdivfree}
\sum_{k\in\cE_d}\wh B_{k,l}(p)
=
\sum_{l\in\cE_d}\wh B_{k,l}(p)
=
0,
\end{align}
and, from \eqref{Chatvector} and \eqref{Chatdivfree} again by direct computations we obtain
\begin{align}
\label{equiv0}
\sum_{k,l\in\cE_d}
(1-e^{-\sqrt{-1}p\cdot k})(1-e^{\sqrt{-1}p\cdot l})\wh B_{k,l}(p)\equiv0.
\end{align}
At $p=0$ we apply $\partial^2/\partial p_i\partial p_j$ to \eqref{equiv0} and get
\begin{align}
\label{C(0)=0}
\wh C_{i,j}(0)
=
\sum_{k,l\in\cE_d}
k_il_j \wh B_{k,l}(0)=0,
\qquad
i,j=1,\dots,d.
\end{align}
Since $\wh C_{j,i}(p) = \wh C_{i,j}(-p) = \overline{\wh C_{i,j}(p)} $
and $p\mapsto \wh C(p)$ is assumed to be twice continuously differentiable at $p=0$, from \eqref{C(0)=0} it follows that
\begin{align}
\notag
\wh C(p)=\Ordo(\abs{p}^2),
\qquad
\text{ as }
\abs{p}\to0,
\end{align}
which implies \eqref{H-1cond_2}.
\end{proof}

\noindent
In particular it follows that sufficiently fast decay of correlations of the divergence-free drift field $V(x)$ implies the ${\cH_{-1}}$-condition \eqref{H-1cond_2}. Note that the divergence-free condition \eqref{divfree} is crucial in this argument.

\section{Historical remarks around Theorem \ref{thm:main} and examples}
\label{s:Historical remarks and examples}

There exist a fair number of important earlier results to which we
should compare Theorem \ref{thm:main}.

\subsubsection*{Remarks}

\begin{enumerate} [(1)]

\item
In Kozlov (1985)  \cite{kozlov_85}, Theorem II.3.3 claims the same result under the supplementary restrictive condition that the drift-field $V(x)$ be \emph{finitely dependent}. However, as pointed out in Komorowski, Olla (2002)  \cite{komorowski_olla_02}, the proof seems to be incomplete there. Also, the condition of finite dependence of the drift field is a very serious restriction.

\item
In Komorowski, Olla (2003)  \cite{komorowski_olla_03a}, Theorem 2.2, essentially the same result is announced as above. However, as noted in section 3.6 of Komorowski et al.\ (2012)  \cite{komorowski_landim_olla_12} this proof is yet again incomplete.

\item
To our knowledge the best fully proved result is Theorem 3.6 of Komorowski et al.\ (2012)  \cite{komorowski_landim_olla_12} where the same result is proved under the condition that the stream tensor field ${H}(x)$ of Proposition \ref{prop:helmholtz} be stationary and in $\cL^{\max\{2+\delta, d\}}$, $\delta>0$, rather than $\cL^2$. Note that the conditions of our theorem only request that the tensor field ${H}$ be square integrable. The proof of Theorem 3.6 in Komorowski et al.\ (2012)  \cite{komorowski_landim_olla_12} is very technical, see sections 3.4 and 3.5 of the monograph.

\item
The special case when the tensor field ${H}$ is actually in $\cL^{\infty}$ is fundamentally simpler. In this case the so-called \emph{strong sector condition} of Varadhan (1996)  \cite{varadhan_96} applies directly. This was noticed in Komorowski, Olla (2003)  \cite{komorowski_olla_03a}. See also section 3.3 of Komorowski et al.\ (2012)  \cite{komorowski_landim_olla_12}, and the examples and comments at the end of this section.

\item
Regarding the continuous space-time diffusion model in divergence-free drift field \eqref{sde}, \eqref{divfree in continuous space}:
It is a fact that, similarly to the $\Z^d$ lattice case, under
minimally restrictive regularity conditions, a stationary and square
integrable divergence-free drift field on $\R^d$ can be written as the
curl of an antisymmetric stream tensor field with stationary
increments ${H}:\R^d\to\R^{d\times d}$. See Proposition 11.1 of
Komorowski et al.\ (2012)  \cite{komorowski_landim_olla_12}, which is the continuous-space
analogue of Proposition \ref{prop:helmholtz}. The
${\cH_{-1}}$-condition \eqref{H-1cond_1}/\eqref{H-1cond_2}/\eqref{H-1cond_3} corresponds to the fact that the
stream tensor ${H}$ is stationary (not just of stationary
increments) and square integrable. The case of bounded ${H}$ was
first considered in Papanicolaou, Varadhan (1981)  \cite{papanicolaou_varadhan_81} --- and this paper
seems to be historically the first instant where the problem of
diffusion in stationary divergence-free drift field was
considered. Homogenization and central limit theorem for the diffusion
\eqref{sde}, \eqref{divfree in continuous space} in \emph{bounded} stream field,
${H}\in\cL^{\infty}$, was first proven in
Papanicolaou, Varadhan (1981)  \cite{papanicolaou_varadhan_81} and Osada (1983)  \cite{osada_83}. The strongest
result in the continuous space-time setup is due to Oelschl\"ager (1988)
 \cite{oelschlager_88} where homogenization and CLT for the
displacement is proved for square-integrable stationary stream tensor
field, ${H}\in\cL^2$.  Oelschl\"ager's proof consists in
truncating the stream tensor and bounding the error. Attempts to apply this
technique in discrete settings run into enormous technical difficulties, see chapter 3 of Komorowski et al.\ (2012)  \cite{komorowski_landim_olla_12}.  If the stream tensor field is stationary \emph{Gaussian} then --- as noted in Komorowski, Olla (2003)  \cite{komorowski_olla_03b} --- the \emph{graded sector condition} of Sethuraman et al.\ (2000)  \cite{sethuraman_varadhan_yau_00} applies. See also chapters 10 and 11 of Komorowski et al.\ (2012) \cite{komorowski_landim_olla_12} for all existing results on the diffusion model \eqref{sde}, \eqref{divfree in continuous space}.

\item
So, our proof of Theorem \ref{thm:main} fills the gap between the restrictive condition ${H}\in\cL^{\max\{2+\delta, d\}}$ of Theorem 3.6 in Komorowski et al.\ (2012)  \cite{komorowski_landim_olla_12} and the minimal restriction ${H}\in\cL^{2}$. The content of our Theorem \ref{thm:main} is the discrete $\Z^d$-counterpart of Theorem 1 in Oelschl\"ager (1988)  \cite{oelschlager_88}. We also stress that our proof is conceptually and technically much simpler that of Theorem 3.6 in Komorowski et al.\ (2012)  \cite{komorowski_landim_olla_12} or Theorem 1 in Oelschl\"ager (1988)  \cite{oelschlager_88}. The continuous space-time diffusion model --- under the same regularity conditions as those of Oelschl\"ager (1988)  \cite{oelschlager_88} can be treated in a very similar way reproducing this way Theorem 1 of Oelschl\"ager (1988)  \cite{oelschlager_88} in a conceptually and technically simpler way. In order to keep this paper relatively short and transparent, those details will be presented elsewhere.

\item
There exist results on \emph{super-diffusive} behaviour of the random walk/diffusion model \eqref{jumpprob}/\eqref{sde} in divergence free drift field, when the ${\cH_{-1}}$-condition \eqref{H-1cond_2} fails to hold. In Komorowski, Olla (2002)  \cite{komorowski_olla_02} and  T\'oth, Valk\'o (2012) \cite{toth_valko_12} the diffusion model \eqref{sde}, \eqref{divfree in continuous space} is considered when the drift field $\Phi$ is Gaussian and the stream tensor field ${H}$ is \emph{genuinely} delocalized: of stationary increment but not stationary. Super-diffusive bounds are proved.

\end{enumerate}

\subsubsection*{Examples}

\newlength{\tempindent}
\setlength{\tempindent}{\parindent}
\begin{enumerate} [leftmargin=0cm,itemindent=0.7cm,labelwidth=\itemindent,labelsep=0cm,align=left,label=(\arabic*)]
\setlength{\parskip}{0cm}\setlength{\parindent}{\tempindent}
\item
\emph{Stationary and bounded stream field:}
When there exists a \emph{bounded} tensor valued variable ${h}:\Omega\to\R^{\cE_{d}\times\cE_{d}}$ with the symmetries \eqref{thetatensor} and such that \eqref{v=curltheta}  holds we define the multiplication operators $N_{k,l}$, $k,l\in\cE_{d}$:
\begin{align}
\label{thetamultipl}
N_{k,l} f(\omega):= {h}_{k,l}(\omega) f(\omega).
\end{align}
These are bounded selfadjoint  operators and they inherit the symmetries of ${h}$:
\begin{align}
\label{Nsymmetries}
\begin{gathered}
N_{l,k}=T_kN_{-k,l}T_{-k}=T_lN_{k,-l}T_{-l}=-N_{k,l},
\\
\sum_{l\in\cE_{d}}N_{k,l}=M_k.
\end{gathered}
\end{align}
As an alternative to \eqref{S and A}, using \eqref{Nsymmetries}, the skew-self-adjoint part of the infinitesimal generator is expressed as
\begin{align}
\label{Aalt}
A=\sum_{k,l\in\cE_{d}} \nabla_{-k}N_{k,l}\nabla_l.
\end{align}
In Komorowski, Olla (2003)  \cite{komorowski_olla_03a} and Komorowski et al.\ (2012)  \cite{komorowski_landim_olla_12} this form of the operator $A$ is used. The operators $N_{k,l}$ are bounded and so is the operator
\begin{align}
\label{b}
B
:=
\abs{\Delta}^{-1/2} A \abs{\Delta}^{-1/2}
=
\sum_{k,l\in\cE_{d}} \Gamma_{-k}N_{k,l}\Gamma_l
\end{align}
which plays a key r\^ole in our forthcoming proof. Due to boundedness of $B$ the \emph{strong sector condition} is valid in these cases and the central limit theorem for the displacement readily follows. See Komorowski, Olla (2003)  \cite{komorowski_olla_03a} and section 3.3 of Komorowski et al.\ (2012) \cite{komorowski_landim_olla_12}.

We give a finitely dependent (actually 1-dependent) constructive example. This construction is in arbitrary dimension --- a corresponding 2-d construction is given in Kozlov (1985)  \cite{kozlov_85}.
Let $x\mapsto ({H}_{e_i,e_j}(x))_{e_i,e_j\in\cE_{d,+}}$ be i.i.d. random vectors such that
\begin{align}
\notag
&
{H}_{e_j,e_i}(x)=-{H}_{e_i,e_j}(x)
&&
\norm{{H}_{e_i,e_j}(x)}_{\infty} < \frac1{2d}.
\end{align}
Define
\begin{align}
\notag
\begin{aligned}
&
{H}_{e_i,-e_j}(x)=-{H}_{e_i,e_j}(x-e_j),
&&
{H}_{-e_i,e_j}(x)=-{H}_{e_i,e_j}(x-e_i),
\\[8pt]
&
{H}_{-e_i,-e_j}(x)={H}_{e_i,e_j}(x-e_i-e_j).
\end{aligned}
\end{align}
The so-called \emph{cyclic walks} analysed in Komorowski, Olla (2003)  \cite{komorowski_olla_03a} and in section 3.3 of Komorowski et al.\ (2012) \cite{komorowski_landim_olla_12} are  also of this nature. As shown in those places, in these cases the \emph{strong sector condition} of Varadhan (1996) \cite{varadhan_96} applies and the central limit theorem for the displacement is obtained.

When the tensor variables ${h}:\Omega\to\R^{\cE_{d}\times\cE_{d}}$ are in $\cL^2\setminus \cL^\infty$ the multiplication operators $N_{k,l}$ defined in \eqref{thetamultipl} are \emph{unbounded}, the representation \eqref{Aalt} of the skew-self-adjoint part of the infinitesimal generator and the operator $B$ defined in \eqref{b} become just \emph{formal}. Nevertheless,  Theorem 1 in Oelschl\"ager (1988)  \cite{oelschlager_88} and theorem 3.6 in Komorowski et al.\ (2012)  \cite{komorowski_landim_olla_12} are proved by controlling approximations of $h_{k,l}$ and the unbounded operators $N_{k,l}$ by truncations at high levels.

\item
\emph{Stationary, square integrable but unbounded stream field:}
We let, in arbitrary dimension $d$, $\Psi:\Z^d+(1/2,\dots, 1/2)\to\Z$ be a stationary, scalar, Lipschitz field with Lipschitz constant 1. As shown in Peled (2010)  \cite{peled_10}, such fields exist in sufficiently high dimension. Define ${H}:\Z^d\to\R^{\cE_{2}\times\cE_{2}}$ by
\begin{align}
\notag
&
{H}_{e_i,e_j}(x)
:=
\frac 1{d}\Psi(x+(e_i+e_j)/2),
&&
\forall x\in\Z^d,
\ \
\forall i<j,
\end{align}
and extend to $\left({H}_{k,l}(x)\right)_{k,l\in\cE_{d}}$ by the symmetries \eqref{Thetasymms}. The tensor field ${H}:\Z^d\to\R^{\cE_{2}\times\cE_{2}}$ defined this way will be stationary and $\cL^2$, but not necessary in $\cL^\infty$ --- the uniform graph homomorphism of Peled (2010)  \cite{peled_10}, for example, is not bounded. Nevertheless, $V$ is bounded by 1, as it should, since $|H_{k,l}(x)+H_{-k,l}(x)|=|H_{k,l}(x)-H_{k,l}(x-k)|\le \frac 1d$ and $V$ is a sum of $d$ such terms.

\item
\emph{Randomly oriented  Manhattan lattice:}
Let $u_i:\Z^{d-1}\to \{-1,+1\}$, $i=1,\dots, d$,  be translation invariant and ergodic, zero mean random fields, which are independent between them. Denote their covariances
\begin{align}
\notag
c_i(y)
&
:=
\expect{u_i(0), u_i(y)},
&&
y\in\Z^{d-1},
\\
\notag
\hat c_i(p)
&
:=
\sum_{y\in\Z^{d-1}} e^{\sqrt{-1}p\cdot y} c_i(y),
&&
p\in[-\pi,\pi)^{d-1}.
\end{align}
Define now the lattice vector field
\begin{align}
\notag
V_{\pm e_i}(x)
:=
\pm  u_i(x_1, \dots, x_{i-1}, \cancel{x_i}, x_{i+1}, \dots, x_d).
\end{align}
Then the random vector field $V$ will satisfy  all conditions in \eqref{lifted conditions} and $t\mapsto X(t)$ will actually be a random walk on the lattice $\Z^d$ whose line-paths parallel to the coordinate axes are randomly oriented in a shift-invariant and ergodic way. This oriented graph is called the \emph{randomly oriented Manhattan lattice}. The covariances $C$ and $\wh C$ defined in \eqref{covariance matrix of Phi}, respectively, \eqref{Fourier transform of covariance matrix of Phi} will be
\begin{align}
\notag
C_{i,j} (x)
&
=
\delta_{i,j} c_i(x_1, \dots, x_{i-1}, \cancel{x_i}, x_{i+1}, \dots, x_d),
\\
\notag
\wh C_{i,j}(p)
&
=
\delta_{i,j} \delta(p_i) \hat c_i(p_1, \dots, p_{i-1}, \cancel{p_i}, p_{i+1}, \dots, p_d).
\end{align}
The $\cH_{-1}$-condition \eqref{H-1cond_2} is in this case
\begin{align}
\notag
\sum_{i=1}^d
\int_{[-\pi,\pi]^{d-1}}
\wh D(q)^{-1}
\hat c_{i}(q) {\d}q <\infty.
\end{align}
In the particular case when the random variables $u_i(y)$, $i\in\{1,\dots, d\}$, $y\in\Z^{d-1}$, are \emph{independent fair coin-tosses}, $\hat c_i(q)\equiv1$. In this case, for $d=2,3$ the ${\cH_{-1}}$-condition  fails to hold, the tensor field ${H}$ is \emph{genuinely} of stationary increments. Super-diffu\-si\-vi\-ty of the walk $t\mapsto X(t)$ is proved like in Tarr\`es et al.\ (2012)  \cite{tarres_toth_valko_12} (in the $2d$ case) and in T\'oth, Valk\'o (2012)  \cite{toth_valko_12} (in the $3d$ case). In dimensions $d\ge4$ the ${\cH_{-1}}$-condition  \eqref{H-1cond_2} holds and the central limit theorem for the displacement can be proved by use of the \emph{graded sector condition} of Sethuraman et al (2000)  \cite{sethuraman_varadhan_yau_00}. If the variables $u_i(y)$, $y\in\Z^{d-1}$ are not independent, the graded sector condition doesn't apply. In this case, however, the central limit theorem follows from our Theorem \ref{thm:main}.

\end{enumerate}

\medskip
\noindent

\section{Relaxed sector condition}
\label{s:Relaxed sector condition}

In  this section we recall from Horv\'ath et al.\ (2012) \cite{horvath_toth_veto_12} the \emph{relaxed sector condition}. This is  a functional analytic condition on the operators $S$ and $A$ which ensures that the efficient martingale approximation \`a la Kipnis-Varadhan of sums of the type of $Z(t)$ in \eqref{Y and Z def} exists. Since in the present case the infinitesimal generator $L=S+A$ is \emph{bounded} we recall the result of Horv\'ath et al.\ (2012) \cite{horvath_toth_veto_12} in a slightly restricted form: we do not have to  worry now about the domains of the operators $S$ and $A$.  This section is fairly abstract.

\subsection{Kipnis-Varadhan theory}
\label{ss:Kipnis-Varadhan theory}

Let $(\Omega, \cF, \pi)$ be a probability space: the state space of a \emph{stationary and ergodic} pure jump Markov process $t\mapsto\eta(t)$ with bounded jump rates. We put ourselves in the Hilbert space $\cL^2(\Omega, \pi)$. Denote the \emph{infinitesimal generator} of the semigroup of the process by $L$. Since the process $\eta(t)$ has bounded jump rates the infinitesimal generator $L$ is a bounded operator.  We denote the \emph{self-adjoint} and \emph{skew-self-adjoint} parts of the generator $L$ by
\begin{align}
\notag
S:=-\frac12(L+L^*)\ge0
\qquad
A:=\frac12(L-L^*).
\end{align}
We assume that $S$ is itself ergodic i.e.
\begin{align}
\notag
\Ker(S)=\{c1\!\!1 : c\in\C\},
\end{align}
and restrict ourselves to the subspace of codimension 1, orthogonal to the constant functions:
\begin{align}
\notag
\cH:=\{f\in\cL^2(\Omega,\pi): \sprod{1\!\!1}{f}=0\}.
\end{align}
In the sequel the operators $(\lambda I+ S)^{\pm1/2}$, $\lambda\ge0$, will play an important r\^ole. These are defined in terms of the spectral theorem applied to the self-adjoint and positive operator $S$.  The \emph{unbounded} operator $S^{-1/2}$ is self-adjoint on its domain
\begin{align}
\notag
\Dom(S^{-1/2})
=
\Ran(S^{1/2})
=
\{f\in\cH:
\norm{S^{-1/2}f}^2
:=
\lim_{\lambda\to0}\norm{(\lambda I + S)^{-1/2}f}^2
<
\infty
\}.
\end{align}
Let $f\in\cH$. We ask about CLT/invariance principle for the rescaled process
\begin{align}
\label{rescaledintegral}
Y_N(t):=\frac{1}{\sqrt{N}}\int_0^{Nt} f(\eta(s)) {\d} s
\end{align}
as $N\to\infty$.

We denote by $R_\lambda$ the \emph{resolvent} of the semigroup $s\mapsto e^{sL}$:
\begin{align}
\notag
R_\lambda
:=
\int_0^\infty e^{-\lambda s} e^{sL} {\d} s
=
\big(\lambda I-L\big)^{-1}, \qquad \lambda>0,
\end{align}
and given $f\in\cH$, we will use the notation
\begin{align}
\notag
u_\lambda:=R_\lambda f.
\end{align}

The following theorem is direct extension to general non-reversible setup of the Kipnis-Varadhan Theorem from  Kipnis, Varadhan (1986)  \cite{kipnis_varadhan_86}. It  yields the \emph{efficient martingale approximation} of the additive functional \eqref{rescaledintegral}.  See T\'oth (1986) \cite{toth_86},  or the surveys Olla (2001)  \cite{olla_01} and Komorowski et al.\ (2012)  \cite{komorowski_landim_olla_12}.

\begin{theoremkv*}
\label{thm:kv}
With the notation and assumptions as before, if for some $v\in\cH$ the following two limits hold in $\cH$:
\begin{align}
\label{conditionA}
\lim_{\lambda\to0}
\lambda^{1/2} u_\lambda=0,
\hskip3cm
\lim_{\lambda\to0} S^{1/2} u_\lambda=v,
\end{align}
then
\begin{align}
\notag
\sigma^2
:=
2\lim_{\lambda\to0}\sprod{u_\lambda}{f}
=
2\norm{v}^2
\in
[0,\infty),
\end{align}
exists, and there also exists a zero mean, $\cL^2$-martingale $M(t)$ adapted to the filtration of the Markov process $\eta(t)$, with stationary and ergodic increments and variance
\begin{align}
\notag
\expect{M(t)^2}=\sigma^2t,
\end{align}
such that
\begin{align}
\notag
\lim_{N\to\infty} N^{-1}
\expect{\abs{Y_N(t) - \frac{M(Nt)}{\sqrt{N}}}^2} =0.
\end{align}
\end{theoremkv*}

\begin{corollarykv*}
With the same setup and notation,
for any $m\in\N$, $t_1,\dots,t_m\in\R_+$ and $F:\R^{m}\to\R$ continuous and bounded
\begin{align}
\notag
\lim_{N\to\infty}
\int_\Omega
\abs{
\expectom{F(Y_N(t_1), \dots, Y_N(t_m))}
-
\expect{F(W(t_1), \dots, W(t_m))}
}
{\d}\pi(\omega)=0,
\end{align}
where $t\mapsto W(t)\in\R$ is a 1-dimensional Brownian motion with variance $\expect{W(t)^2}=\sigma^2 t$.
\end{corollarykv*}

\subsection{Relaxed sector condition}
\label{ss:Relaxed sector condition}

Let, for $\lambda>0$,
\begin{align}
\label{Blambda_def}
B_\lambda:=(\lambda I + S)^{-1/2} A (\lambda I + S)^{-1/2}.
\end{align}
These are bounded and skew-self-adjoint.

The main result of Horv\'ath et al.\ (2012) \cite{horvath_toth_veto_12}, restricted to the case when the operators $S$ and $A$ are bounded, is the following:

\begin{theoremrsc*}
\label{thm:rsc}
Assume that there exist a dense subspace $\cC\subseteq \cH$ and an operator $B:\cC\to \cH$ which is essentially skew-self-adjoint on the core $\cC$ and such that for any vector $\varphi\in \cC$
\begin{align}
\label{Blambdalimit}
\lim_{\lambda\to 0}\norm{B_\lambda\varphi-B\varphi}=0.
\end{align}
Then, the  ${\cH_{-1}}$-condition
\begin{align}
\notag
f\in\Dom(S^{-1/2})
\end{align}
implies    \eqref{conditionA}, and thus the martingale approximation and CLT of Theorem KV follow.
\end{theoremrsc*}

The unbounded operator $B$ is  \emph{formally} identified as
\begin{align}
\notag
B:=S^{-1/2} A S^{-1/2}.
\end{align}

RSC refers to \emph{relaxed sector condition}: indeed, as shown in
Horv\'ath et al.\ (2012) \cite{horvath_toth_veto_12} this theorem contains the \emph{strong
  sector condition} of Varadhan (1996)  \cite{varadhan_96} and the \emph{graded sector
  condition} of Sethuraman et al.\ (2000)  \cite{sethuraman_varadhan_yau_00} as special
cases. For comments on history, content and variants of Theorem KV we
refer the reader to the  monograph
Komorowski et al.\ (2012)  \cite{komorowski_landim_olla_12}. For the full proof and direct consequences of Theorem RSC see Horv\'ath et al.\ (2012) \cite{horvath_toth_veto_12}.

\section{The operator \texorpdfstring{$B=S^{-1/2}AS^{-1/2}$}{B} and
  proof of Theorem \ref{thm:main}}
\label{s:The operator B and proof of Theorem 1}

Let
\begin{align}
\notag
\cC
:=
\Dom\nsqd
=
\Ran\sqd,
\end{align}
and recall from \eqref{Gamma_k} and \eqref{vmultipl} the definition of the operators $\Gamma_l$ and $M_l$, $l\in\cE_d$. Define the \emph{unbounded} operator $B:\cC\to\cH$
\begin{align}
\notag
B
:=
-
\sum_{l\in\cE_{d}}\Gamma_{-l} M_l \nsqd.
\end{align}
First we verify \eqref{Blambdalimit}, i.e.\ that $B_\lambda\to B$ \emph{pointwise} on the core $\cC$, where the bounded operator  $B_\lambda$ is expressed by inserting the explicit form of $S$ and $A$ \eqref{S and A} into the definition of $B_\lambda$ \eqref{Blambda_def}:
\begin{align}
\notag
B_\lambda = -\sum_{l\in\cE^d}(\lambda I - \Delta)^{-1/2}\nabla_{-l}M_l(\lambda I - \Delta)^{-1/2}.
\end{align}
From the spectral theorem for the commutative $C^*$-algebra generated by the shift operators $T_{e_i}$, $i=1,\dots,d$, (see e.g.\ Arveson (1976) \cite{arveson_76}) we obtain that $\|(\lambda I-\Delta)^{-1/2}\nabla_l\|\le 1$, $\|(\lambda I-\Delta)^{-1/2}\abs{\Delta}^{1/2}\|\le 1$, and, moreover, for any $\varphi\in\cH$
\begin{align}
\notag
&
(\lambda I - \Delta)^{-1/2} \nabla_l \varphi
\to
\Gamma_l \varphi
&&
(\lambda I - \Delta)^{-1/2} \abs{\Delta}^{1/2} \varphi
\to
\varphi.
\end{align}
Hence \eqref{Blambdalimit} follows readily for any $\varphi\in\cC$.

With \eqref{Blambdalimit} established, we need to show that $B$ is essentially skew-self-adjoint on $\cC$. We start with a light lemma.

\begin{lemma}
\label{lem:B^*}
(i)
$B:\cC\to\cH$ is skew-Hermitian.
\\
(ii)
The full domain of $B^*$ is
\begin{align}
\label{domain of B^*}
{\cC}^*
=
\{f\in\cH: \sum_{l\in\cE_{d}} M_l \Gamma_l f \in\cC\}
\end{align}
and $B^*$ acts on $\cC^*$ by
\begin{align}
\label{B^*}
B^*
:=
-
\nsqd \sum_{l\in\cE_{d}} M_l \Gamma_l.
\end{align}
\end{lemma}

\bigskip
\noindent
{\bf Remark:}
It is of crucial importance here that $\cC^*$ in \eqref{domain of B^*}
is the \emph{full} domain of the adjoint operator $B^*$, i.e.~the subspace
of all $f$ such that the linear functional $g\mapsto \sprod{f}{Bg}$ is bounded on
$\cC$. It will not be enough for our purposes just to show that
$\cC^*$ is some core of definition.

\bigskip

\begin{proof}
(i)
Let $f,g\in\cC$. Then, due to \eqref{Mnablacommute}
\begin{align}
\notag
\sprod{f}{Bg}
&
=
-
\sum_{l\in\cE_{d}} \sprod{\nsqd f}{\nabla_{-l} M_l \nsqd g}
\\
\notag
&
=
\phantom{-}
\sum_{l\in\cE_{d}} \sprod{\nabla_{-l} M_l \nsqd f}{\nsqd g}
=
-\sprod{Bf}{g},
\end{align}
(ii)
Next,
\begin{align}
\notag
\Dom(B^*)
&=
\Big\{
f\in\cH:
(\exists c(f)<\infty)
(\forall g\in\cC):
\Big|\Big\langle f, \sum_{l\in\cE_d} \Gamma_{-l}M_l \nsqd g\Big\rangle\Big|
\le
c(f)\norm{g}
\Big\}
\\
\notag
&=
\Big\{
f\in\cH:
(\exists c(f)<\infty)
(\forall g\in\cC):
\Big|\Big\langle\sum_{l\in\cE_d} M_l \Gamma_{l} f,\nsqd g\Big\rangle\Big|
\le
c(f)\norm{g}
\Big\}
\\
\notag
&
=
\Big\{
f\in\cH:
\sum_{l\in\cE_d} M_l \Gamma_{l} f \in \cC
\Big\},
\end{align}
as claimed. In the last step we used the fact that $\cC$ is the
\emph{full domain} of the self-adjoint operator $\nsqd$. The action \eqref{B^*} of $B^*$ follows from straightforward manipulations.
\end{proof}

Note that Lemma \ref{lem:B^*} in particular implies that $\cC\subseteq\cC^*$, the operator $B:\cC\to\cH$ is closable and $B^*:\cC^*\to\cH$ is extension of $-\overline{B}$. We actually ought to prove that
\begin{align}
\notag
B^*=-\overline{B}.
\end{align}
In order to do this we still need to check von Neumann's criterion (see
e.g.\ Theorem VIII.3 of Reed, Simon (1975)  \cite{reed_simon_vol1_vol2_75}):
\begin{align}
\label{vonNeumann's criterion}
\Ker(B^* \pm I)=\{0\},
\end{align}
or, equivalently, prove that the equations
\begin{align}
\label{vonNeumannplus}
\sum_{l\in\cE_{d}}M_l \Gamma_l \mu + \sqd \mu = 0
\\
\label{vonNeumannminus}
\sum_{l\in\cE_{d}}M_l \Gamma_l \mu - \sqd \mu = 0
\end{align}
admit only the trivial solution $\mu=0$. We will prove this for \eqref{vonNeumannplus}. The other case is done very similarly.

Note that assuming $\mu\in\cC$ the problem becomes fully trivial. Indeed, inserting $\mu=\abs{\Delta}^{1/2}\chi$ in \eqref{vonNeumannplus} the equation is rewritten as
\begin{align}
\notag
\sum_{l\in\cE_d} M_l \nabla_l \chi - \Delta \chi =0,
\end{align}
which admits only the trivial solution $\chi=0$, due to  \eqref{kerDelta} and \eqref{Mnablacommute}. The point is that $\mu$ is not necessarily in $\cC$ so $|\Delta|^{-1/2}\mu$ is not necessarily well defined as an element of $\cH$. Nevertheless, we are able to define a scalar random field $\Psi:\Omega\times\Z^d \to \R$ of stationary increments (rather than stationary) which can be thought of as the lifting of $|\Delta|^{-1/2}\mu$ to the lattice $\Z^d$.

Let, therefore $\mu$ be a putative solution for \eqref{vonNeumannplus} and define, for each $k\in\cE_d$,
\begin{align}
\label{udef}
u_k:=\Gamma_k \mu.
\end{align}
These are vector components and they also satisfy the \emph{gradient condition}: $\forall \ k,l\in\cE_d$
\begin{align}
\notag
&
u_k+T_ku_{-k}=0,
&&
u_k+T_ku_l=u_l+T_lu_k.
\end{align}
Note also that
\begin{align}
\notag
\sum_{l\in\cE_{d}} u_l = \sqd \mu.
\end{align}
The eigenvalue equation \eqref{vonNeumannplus} becomes
\begin{align}
\label{vonNeumann3}
\sum_{l\in\cE_{d}} v_l u_l + \sum_{l\in\cE_{d}}u_l=0.
\end{align}
We lift this equation to $\Z^d$, by \eqref{vlift}, and similarly defining the lattice vector field $U:\Omega\times\Z^d\to\R^d$ as
\begin{align}
\notag
U_k(\omega,x):=u_k(\tau_x\omega).
\end{align}
Then $U$ is the $\Z^d$-gradient of a scalar field $\Psi:\Omega\times\Z^d\to\R$, determined uniquely by
\begin{align}
\label{Psidef}
&
\Psi(\omega,0)=0,
&&
\Psi(\omega, x+k) -\Psi(\omega, x) = U_k(\omega, x).
\end{align}
Note that, as promised, the scalar field $\Psi$ has stationary increments (or, in the language of ergodic theory: it is a cocycle).
\begin{align}
\label{Psistatincr}
\Psi(\omega, x) - \Psi(\omega, y) = \Psi(\tau_x\omega, y-x) - \Psi(\tau_x\omega, 0).
\end{align}
The lifted version of the equation \eqref{vonNeumann3} is
\begin{align}
\label{harmonic1}
\sum_{l\in\cE_{d}}(\Psi(\omega,x+l)-\Psi(\omega,x))
+
\sum_{l\in\cE_{d}}V_l(\omega,x)(\Psi(\omega,x+l)-\Psi(\omega,x))
=0,
\end{align}
or, shortly,
\begin{align}
\notag
\lap \Psi +  V\cdot \grad \Psi = 0,
\end{align}
where here $\lap$ and $\grad$ denote the $\Z^d$-Laplacian and $\Z^d$-gradient. In view of \eqref{jumpprob} equation \eqref{harmonic1} means exactly that given $\omega\in\Omega$ fixed (that is: in the quenched setup) the field $\Psi:\Z^d\to\R$ is \emph{harmonic} for the random walk $X(t)$. Thus the process
\begin{align}
\label{marti}
t\mapsto R(t):=\Psi(X(t))
\end{align}
is a martingale (with $R(0)=0$) in the quenched filtration $\sigma\left( \omega,X(s)_{0\le s\le t} \right)_{t\ge0}$.
On the other hand, from stationarity and ergodicity of the environment process $t\mapsto \eta_t$ and \eqref{Psistatincr} it follows that the process $t\mapsto R(t)$ has stationary and ergodic increments with respect to the annealed measure $\probab{\cdot}:= \int_\Omega \probabom{\cdot} d\pi(\omega)$. Indeed, using \eqref{Psistatincr}, a straightforward computation shows that
\begin{align}
\notag
\expectom{F(R(t_0+\cdot)-R(t_0))}
=
{\mathbf{E}_{\eta_{t_0}}\left(F(R(\cdot)\right)},
\end{align}
where $F(R(\cdot))$ is an arbitrary bounded and measurable functional of the process $t\mapsto R(t)$, $t\ge0$. Hence, by stationarity of the environment $t\mapsto\eta_t$, the claim follows.

Hence it follows that the process $t\mapsto R(t)$ is a martingale (with $R(0)=0$) with stationary and ergodic increments, in its own annealed filtration $\sigma\left(R(s)_{0\le s\le t} \right)_{t\ge0}$.

\begin{lemma}
\label{lem:variance of R}
Let $\mu$ be solution of the equation \eqref{vonNeumannplus}, $\Psi$ the harmonic field constructed in \eqref{Psidef} and $R(t)$ the martingale defined in \eqref{marti}. Then
\begin{align}
\label{martingale variance}
\expect{R(t)^2} = 2\norm{\mu}^2 t.
\end{align}
\end{lemma}

\begin{proof}
Since $t\mapsto R(t)$ is a martingale with stationary increments (with respect to the annealed measure $\probab{\cdot}$), we automatically have $\expect{R(t)^2}=\varrho^2 t$ with some $\varrho\ge0$. We now compute $\varrho$.
\begin{align}
\notag
\varrho^2
&
:=
\lim_{t\to0}
\frac{\expect{R(t)^2}}{t}
\stackrel{(1)}{=}
\lim_{t\to0} \int_\Omega
\frac{\expectom{\Psi(\omega, X(t))^2}}{t}
\d\pi(\omega)
\\
\notag
&
\stackrel{(2)}{=}
\int_\Omega \lim_{t\to0}
\frac{\expectom{\Psi(\omega, X(t))^2}}{t}
\d\pi(\omega)
\stackrel{(3)}{=}
\sum_{l\in\cE_{d}}\int_{\Omega}
\left(1 + v_l(\omega)\right) \abs{u_l(\omega)}^2
{\d}\pi(\omega)
\\
\notag
&
\stackrel{(4)}{=}
\sum_{l\in\cE_{d}}\int_{\Omega}
\abs{u_l(\omega)}^2
{\d}\pi(\omega)
\stackrel{(5)}{=}
\sum_{l\in\cE_{d}}
\norm{\Gamma_l\mu}^2
\stackrel{(6)}{=}
2 \norm{\mu}^2.
\end{align}
Step (1) is annealed averaging.
Step (2) is easily justified by dominated convergence.
Step (3) drops out from explicit computation of the conditional variance of one jump.
In step (4)  we used that $v_{-l}(\omega)\abs{u_{-l}(\omega)}^2 = -
v_{l}(\tau_{-l}\omega)\abs{u_{l}(\tau_{-l}\omega)}^2$ and translation
invariance of the measure $\pi$ on $\Omega$.
In step (5) we use the definition \eqref{udef} of $u_l$.
Finally, in the last step (6) we used the third
identity of \eqref{Gammaadj}.
\end{proof}

\begin{proposition}
\label{prop:Neumann}
The unique solution of \eqref{vonNeumannplus}/\eqref{vonNeumannminus} is $\mu=0$, and consequently the operator $B$ is essentially skew-self-adjoint on the core $\cC$.
\end{proposition}

\begin{proof}
Let $\mu$ be a solution of the equation \eqref{vonNeumannplus}, $\Psi$ the harmonic field constructed in \eqref{Psidef} and $R(t)$ the martingale defined in \eqref{marti}. From the martingale central limit theorem and \eqref{martingale variance} it follows that
\begin{align}
\label{martingale CLT}
\frac{R(t)}{\sqrt{t}}
\Rightarrow \cN(0,2\norm{\mu}^2),
\qquad
\text{ as } t\to\infty.
\end{align}
On the other hand we are going to prove that
\begin{align}
\label{toprob0}
\frac{R(t)}{\sqrt{t}}
\toprob0,
\qquad
\text{ as }t\to\infty.
\end{align}
Jointly, \eqref{martingale CLT} and \eqref{toprob0} clearly imply $\mu=0$, as claimed in the proposition.

The proof will combine
\\
(A)
the (sub)diffusive behaviour of the displacement
\begin{align}
\label{X is (sub)diffusive}
\varlimsup_{T\to\infty}T^{-1}\expect{X(T)^2}<\infty,
\end{align}
which follows from the ${\cH_{-1}}$-condition,  see \eqref{diffusive bounds}; and
\\
(B)
the fact that the scalar field $x\mapsto \Psi(x)$ having zero mean and stationary increments, cf. \eqref{Psistatincr}, increases \emph{sublinearly} with $\abs{x}$. The sublinearity is the issue here. Since $\Psi$ has stationary, mean zero increments, due to the individual (pointwise) ergodic theorem, it follows that \emph{on any fixed line} $\Psi$ increases sublinearly almost surely. However, this does not warrant that $\Psi$ increases sublinearly uniformly in $\Z^d$, $d\ge 2$ which is the difficulty we will now tackle.

Let $\delta>0$ and $K<\infty$. Then
\begin{align}
\label{R large}
\probab{\abs{R(t)}>\delta\sqrt{t}}
\le
\probab{\{\abs{R(t)}>\delta\sqrt{t}\} \cap \{\abs{X(t)}\le K\sqrt{t}\}}
+
\probab{\abs{X(t)}> K\sqrt{t}}.
\end{align}
From (sub)diffusivity \eqref{diffusive bounds} and Chebyshev's inequality it follows directly that
\begin{align}
\label{X large}
\lim_{K\to\infty}
\varlimsup_{t\to\infty}
\probab{\abs{X(t)}> K\sqrt{t}}
=0.
\end{align}
We present two proofs of
\begin{align}
\label{R large and X small}
\lim_{t\to\infty}
\probab{\{\abs{R(t)}>\delta\sqrt{t}\} \cap \{\abs{X(t)}\le K\sqrt{t}\}}
=0,
\end{align}
with $K<\infty$ fixed.  One with bare hands, valid in $d=2$ only, and another one valid in any  dimension which relies on a heat kernel (upper) bound from Morris, Peres (2005)  \cite{morris_peres_05}.

\begin{proof}
[Proof of \eqref{R large and X small} in $d=2$, with bare hands]
First note that
\begin{align}
\label{bound R by Psi}
\probab{\{\abs{R(t)}>\delta\sqrt{t}\} \cap \{\abs{X(t)}\le K\sqrt{t}\}}
\le
\probab{ \max_{\abs{x}<K\sqrt{t}} \abs{\Psi(x)}> \delta\sqrt{t}}.
\end{align}
Next, since $\Psi$ is harmonic with respect to the random walk $X(t)$, it obeys the \emph{maximum principle} (this is true for any random walk, no special property of $X$ is used here). Thus
\begin{align}
\label{maximum principle}
\max_{\abs{x}_\infty\le L}\abs{\Psi(x)}
=
\max_{\abs{x}_\infty = L}\abs{\Psi(x)},
\end{align}
where $\abs{x}_\infty:=\max\{\abs{x_1},\abs{x_2}\}$.  By spatial stationarity
\begin{align}
\label{shiftit}
\begin{gathered}
\max_{\abs{x_1}\le L} \abs{\Psi(x_1,-L)-\Psi(0,-L)}
\sim
\max_{\abs{x_1}\le L} \abs{\Psi(x_1,0)}
\sim
\max_{\abs{x_1}\le L} \abs{\Psi(x_1,+L)-\Psi(0,+L)},
\\
\max_{\abs{x_2}\le L} \abs{\Psi(-L,x_2)-\Psi(-L,0)}
\sim
\max_{\abs{x_2}\le L} \abs{\Psi(0,x_2)}
\sim
\max_{\abs{x_2}\le L} \abs{\Psi(+L,x_2)-\Psi(+L,0)},
\end{gathered}
\end{align}
where $\sim$ stands for equality in distribution. Now, note that $\Psi(x_1,0)$ and $\Psi(0, x_2)$ are Birkhoff sums:
\begin{align}
\notag
\Psi(x_1,0)
=
\sum_{j=0}^{x_1-1} u_{e_1}(\tau_{je_1}\omega),
\qquad
\Psi(0, x_2)
=
\sum_{j=0}^{x_2-1} u_{e_2}(\tau_{je_2}\omega),
\end{align}
where $u_{e_1}(\omega)$ and $u_{e_1}(\omega)$ are zero mean and square integrable.
Hence, by the ergodic theorem
\begin{align}
\label{ergthm}
L^{-1}
\max\big\{
\max_{\abs{x_1}\le L} \abs{\Psi(x_1,0)},
\max_{\abs{x_2}\le L} \abs{\Psi(0,x_2)}
\big\}
\toas0,
\qquad
\text{ as }
L\to\infty.
\end{align}
Putting together \eqref{maximum principle}, \eqref{shiftit} and \eqref{ergthm} we readily obtain, for any $\vareps>0$,
\begin{align}
\label{Psi is sublinear}
\lim_{L\to\infty}
\probab{\max_{\abs{x}_\infty\le L}\abs{\Psi(x)}\ge\vareps L}
=0.
\end{align}
Finally, \eqref{R large and X small} follows by applying \eqref{Psi is sublinear} to the right hand side of \eqref{bound R by Psi}.
\end{proof}

\begin{proof}
[Proof of \eqref{R large and X small} in all $d\ge 2$]
We start with the following uniform upper bound on the (quenched) heat kernel of the walk $X(t)$.

\begin{proposition}
\label{prop:heat kernel bound}
There exists a constant $C=C(d)$ (depending only on the dimension $d$) such that for $\pi$-almost all $\omega\in\Omega$ and all $t>0$
\begin{align}
\label{heat kernel bound}
\sup_{x\in\Z^d} \probabom{X(t)=x}\le C t^{-d/2}.
\end{align}
\end{proposition}

\begin{proof}
This bound \eqref{heat kernel bound} follows from Theorem 2 of Morris, Peres (2005)  \cite{morris_peres_05} through Lemma \ref{lem:Morris-Peres}, below, which states essentially the same bound for discrete-time lazy random walks on $\Z^d$.

\begin{lemma}
\label{lem:Morris-Peres}
Let $V:\cE_d\times\Z^d\to[-1,1]$ be a (deterministically given) field satisfying the conditions in \eqref{lifted conditions} and $n\mapsto X_n$ a discrete-time nearest-neighbour, lazy random walk on $\Z^d$ with transition probabilities
\begin{align}
\notag
\condprobab{X_{n+1}=y}{X_n=x}
=
p_{x,y}
:=
\begin{cases}
\frac12
&\text{ if }
y=x,
\\
\frac{1}{4d}(1+V_k(x))
&\text{ if }
y=x+k, \ \  k\in\cE_d,
\\
0
&\text{ if }
\abs{y-x}>1.
\end{cases}
\end{align}
There exists a constant $C=C(d)$ depending only on dimension such that for any $x,y\in\Z^d$
\begin{align}
\label{discrete heat kernel bound}
\condprobab{X_n=y}{X_0=x}\le C n^{-d/2}.
\end{align}
\end{lemma}

\begin{proof}
For $A,B\subset\Z^d$, such that $A\cap B=\emptyset$ let
\begin{align}
\notag
Q(A,B)
:=
\sum_{x\in A, y\in B}
p(x,y).
\end{align}
For $S\subset\Z^d$, $\abs{S}<\infty$
\begin{align}
\notag
Q(S,S^c)
&
=
\sum_{x\in S, y\in S^c} \frac{1}{4d} (1+V_{y-x}(x))
\\
\notag
&
=
\frac{1}{4d} \abs{\partial S}
+
\frac{1}{4d}
\left(
\sum_{x\in S, y\in\Z^d}
V_{y-x}(x)
-
\sum_{x\in S, y\in S}
V_{y-x}(x)
\right)
\\
\label{cut}
&
=
\frac{1}{4d} \abs{\partial S},
\end{align}
where the last equality follows from
\begin{align}
\notag
\sum_{x\in S, y\in\Z^d}
V_{y-x}(x)
&
=
\sum_{x\in S} \sum_{l\in\cE_d} V_l(x) =0,
\\
\notag
\sum_{x\in S, y\in S}
V_{y-x}(x)
&
=
\frac12
\sum_{x\in S, y\in S}
\left(
V_{y-x}(x)
+
V_{x-y}(y)
\right)
=0,
\end{align}
both of which are consequences of \eqref{lifted conditions}. Since the uniform counting measure on $\Z^d$ is stationary to our walk (from divergence-freeness, cf.  \eqref{lifted conditions}), the isoperimetric profile $\Phi(r)$ (in the sense of Morris, Peres (2005)  \cite{morris_peres_05}) is defined by
\begin{align}
\notag
\Phi(r)
:=
\inf_{0<\abs{S}\le r}\frac{Q(S,S^c)}{\abs{S}}.
\end{align}
Theorem 2 of Morris, Peres (2005)  \cite{morris_peres_05} (specified to our setup) states that for any $0<\vareps\le 1$, if
\begin{align}
\label{morris-peres thm2}
n> 1 + 4 \int_{4}^{4/\vareps}\frac{{\d} u}{u\Phi^2(u)}
\end{align}
then, for any $x, y \in \Z_d$
\begin{align}
\notag
\condprobab{X_n=y}{X_0=x} \le \vareps.
\end{align}
From \eqref{cut} and the standard isoperimetric inequality on $\Z^d$ we have
\begin{align}
\label{isope}
C_1 r^{-1/d} \le \Phi(r) \le C_2 r^{-1/d},
\end{align}
with the constants $0<C_1<C_2<\infty$ depending only on the dimension.
Finally, from \eqref{morris-peres thm2} and \eqref{isope} we readily get \eqref{discrete heat kernel bound}.
\end{proof}

In order to obtain \eqref{heat kernel bound} from \eqref{discrete heat kernel bound}, note that $X(t)=X_{\nu(t)}$ where $\nu(t)$ are distributed like Poisson random variables with parameter $t/2$, and are independent of the discrete time walk $X_n$. Thus
\begin{align}
\notag
\probabom{X(t)=x}
&
=
e^{-t/2}\sum_{n=0}^\infty \frac{(t/2)^n}{n!} \probabom{X_n=x}
\\
\notag
&
\le
e^{-t/2}\sum_{n=0}^\infty \frac{(t/2)^n}{n!} Cn^{-d/2}
\\
\notag
&
\le
C t^{-d/2}
\end{align}
This completes the proof of Proposition \ref{prop:heat kernel bound}.
\end{proof}

\subsubsection*{Remarks.}
\setlength{\tempindent}{\parindent}
\begin{enumerate} [leftmargin=0cm,itemindent=0.7cm,labelwidth=\itemindent,labelsep=0cm,align=left,label=(\arabic*)]
\setlength{\parskip}{0cm}\setlength{\parindent}{\tempindent}
\item
The point in Proposition \ref{prop:heat kernel bound} is that it provides uniform upper bound in any (deterministic) environment which satisfies conditions \eqref{lifted conditions}.

\item
A similar statement and proof holds for the more general case of random walk in doubly stochastic environment as defined in \eqref{jumpprob_gen}, \eqref{bistochastic}, with the symmetric part of the jump probabilities satisfying the uniform ellipticity condition \eqref{ellipticity}.

\item
In Lemma \ref{lem:Morris-Peres} the ``amount of laziness'' could be any $\delta\in(0,1)$, with appropriate minor changes in the formulation and proof.

\item
Alternative proofs of Proposition \ref{prop:heat kernel bound} are also valid, using either Nash-Sobolev or Faber-Krahn inequalities, see e.g. Kumagai (2014) \cite{kumagai_14}. These alternative proofs -- which we do not present here -- are more analytic in flavour. Their advantage is robustness: these proofs are also valid in the continuous-space setting of \eqref{sde}, \eqref{divfree in continuous space}.
\end{enumerate}

We now return to the proof of \eqref{R large and X small}.  By Chebyshev's inequality
\begin{align}
\label{cheb}
\probab{\{\abs{R(t)}>\delta\sqrt{t}\} \cap \{\abs{X(t)}\le K\sqrt{t}\}}
\le
\delta^{-2}t^{-1}\expect{\abs{R(t)}^2\ind{\abs{X(t)}\le K\sqrt{t}}}
\end{align}
Since the scalar field $\Psi$ has stationary increments,
cf. \eqref{Psistatincr}, and zero mean, we get from the $\cL^2$ ergodic theorem that
\begin{align}
\notag
\lim_{n\to\infty}n^{-2}\expect{\abs{\Psi(ne_k)}^2} =0\qquad\forall e_k\in\cE_d,
\end{align}
and, consequently,
\begin{align}
\label{Psi sublinear}
\lim_{|x|\to\infty}|x|^{-2}\expect{\abs{\Psi(|x|)}^2} =0.
\end{align}
Applying in turn the heat kernel bound \eqref{heat kernel bound} of Proposition \ref{prop:heat kernel bound} and the limit \eqref{Psi sublinear} on the right hand side of \eqref{cheb} we obtain
\begin{align}
\notag
t^{-1}\expect{\abs{R(t)}^2 \ind{\abs{X(t)}\le K\sqrt{t}}}
\le
Ct^{-d/2-1}\sum_{|x|\le K\sqrt{t}}\expect{|\Psi(x)|^2}
\to
0,
\qquad
\textrm{ as }t\to\infty.
\end{align}
Here the first expectation is both on the random walk $X(t)$ and on the
field $\omega$, while the second is just on the field $\omega$. The point is that with the help of the uniform heat kernel bound of Proposition \ref{prop:heat kernel bound} we can \emph{decouple} the two expectations.

This concludes the proof of \eqref{R large and X small} in arbitrary dimension.
\end{proof}

\noindent
We conclude the proof of the Proposition \ref{prop:Neumann} by noting that from \eqref{R large}, \eqref{X large} and \eqref{R large and X small} we readily get \eqref{toprob0} which, together with \eqref{martingale CLT} implies indeed that $\mu=0$. So \eqref{vonNeumann's criterion} holds. (We showed that $\Ker(B^*+I)=\{0\}$, the proof $\Ker(B^*-I)=\{0\}$ is the same with $-V$ instead of $V$). Thus $B$ is indeed skew-self-adjoint.
\end{proof}

\medskip
\noindent
{\bf Remark.}
In the proof of Proposition \ref{prop:Neumann} we only use the upper bound $\expect{\abs{X(t)}^2}\le Ct$ which is a consequence \eqref{X is (sub)diffusive} of the ${\cH_{-1}}$-condition \eqref{H-1cond_2}, but is not equivalent to it.

\medskip
\noindent
Altogether we have proved the following theorem, slightly more general than Theorem \ref{thm:main}

\begin{theorem}
\label{thm:more general}
Assume that the diffusive upper bound \eqref{X is (sub)diffusive} on the displacement of the random walker holds. Than, for any $f\in\Dom(|\Delta|^{-1/2})$ the efficient martingale approximation of Theorem KV is valid.
\end{theorem}

Finally, Theorem \ref{thm:main} follows from the ${\cH_{-1}}$-condition \eqref{H-1cond_2}, the bounds \eqref{bounds on sigmasquared} and Theorem \ref{thm:more general}.
\qed

\bigskip\bigskip

\noindent
{\bf Acknowledgements:}
The research of BT is partially supported by OTKA (Hungarian National Research Fund) grant K 100473. The research of GK partially supported by the Israel Science Foundation. Both authors acknowledge mobility support by The Leverhulme Trust through the International Network ``Laplacians, Random Walks, Quantum Spin Systems''.

\vskip2cm

\hbox{
\vbox{\hsize=7cm\noindent
{\sc Gady Kozma}
\\
Department of Mathematics
\\
The Weizmann Institute of Science
\\
POB 26, Rehovot, 76100
\\
Israel
\\
email: {\tt gady.kozma@weizmann.ac.il}
}
\hskip2cm
\vbox{\hsize=7cm\noindent
{\sc B\'alint T\'oth}
\\
School of Mathematics
\\
University of Bristol
\\
Bristol, BS8 1TW
\\
United Kingdom
\\
email: {\tt balint.toth@bristol.ac.uk}
}
}

\end{document}